\documentclass{article}

\usepackage[utf8]{inputenc}
\usepackage{graphicx}
\usepackage{multirow}
\usepackage{amssymb}
\usepackage{amsthm}
\usepackage{amsmath}
\usepackage{algorithm, algorithmicx, algpseudocode}
\newtheorem{thm}{Theorem}
\newtheorem{lem}[thm]{Lemma}
\newtheorem{cor}[thm]{Corollary}

\usepackage[top=1in, bottom=1in, left=1in, right=1in]{geometry}

\begin{document}

\title{A stabilized GMRES method for singular and severely ill-conditioned systems of linear equations}

\author{Zeyu LIAO\footnote{Department of Informatics, School of Multidisciplinary Sciences, The Graduate University for Advanced Studies (SOKENDAI), 2-1-2 Hitotsubashi, Chiyoda-ku, Tokyo, 101-8430, Japan,
                          Email: zeyu@nii.ac.jp}, \:         
        Ken HAYAMI\footnote{Professor Emeritus, National Institute of Informatics, and The Graduate University for Advanced Studies (SOKENDAI),  Email: hayami@nii.ac.jp}, \:
Keiichi MORIKUNI\footnote{Faculty of Engineering, Information and Systems, University of Tsukuba, 1-1-1 Tennodai, Tsukuba, Ibaraki, 305-8573, Japan, Email: morikuni@cs.tsukuba.ac.jp}, \, and \, Jun-Feng YIN\footnote{School of Mathematical Science, Tongji University,  Siping Road 1239, Yangpu District, Shanghai, 200092,  China, Email: yinjf@tongji.edu.cn} 
}

\date{}

\maketitle

\begin{abstract}
Consider using the right-preconditioned GMRES (AB-GMRES) for obtaining the minimum-norm solution of inconsistent underdetermined systems of linear equations. Morikuni (Ph.D. thesis, 2013) showed that for some inconsistent and ill-conditioned problems, the iterates may diverge. This is mainly because the Hessenberg matrix in the GMRES method becomes very ill-conditioned so that the backward substitution of the resulting triangular system becomes numerically unstable.
 We propose a stabilized GMRES based on solving the normal equations corresponding to the above triangular system using the standard Cholesky decomposition. This has the effect of shifting upwards the tiny singular values of the Hessenberg matrix which lead to an inaccurate solution.
We analyze why the method works.
Numerical experiments show that the proposed method is robust and efficient, not only for applying AB-GMRES to underdetermined systems, but also for applying GMRES to
severely ill-conditioned range-symmetric systems of linear equations.
\end{abstract}

\vspace{5mm}

\section{Introduction}\label{sec1}
As a motivating instance when the generalized minimal residual (GMRES) method iterates diverge due to severe ill-conditioning, consider obtaining the minimum-norm solution of  the inconsistent least squares problem:
\begin{equation}\label{eq3}
\min_{x\in \mathbb{R}^n}\|x\|_2, \rm{such\ that}\ \it x\in \{\arg \min_{\xi\in\mathbb{R}^n}\|b-A\xi\|_{\rm{2}}\}
\end{equation}
where $A\in \mathbb{R}^{m\times n}$ and $b\notin \rm{\mathcal{R}}(\it{A})\subseteq \mathbb{R}^{m}$. Here, $\rm{\mathcal{R}}$$(A)$ denotes the range space of $A$. Such problems may occur in ill-posed problems where $b$ is given by an observation which contains noise.
The problem
(\ref{eq3}) is equivalent to
\begin{equation}
(A^{\mathsf{T}}A)^2v=A^{\mathsf{T}}b, x=A^{\mathsf{T}}Av,
\end{equation}
and the solution  can be expressed by $x=A^\dag b$, where $A^{\mathsf{T}}$ denotes the transpose of $A$ and $A^\dag$ is the pseudoinverse of $A$. (See e.g. \cite{bjorck1996numerical}.)\par
The standard direct method for solving the least squares problem $(\ref{eq3})$ is to use the QR decomposition. However, when $A$ is large and sparse, iterative methods become necessary. The CGLS \cite{HS} and LSQR \cite{lsqr} are mathmetically equivalent to applying the conjugate gradient (CG) method to the normal equations of the first kind  \begin{equation}\label{eq4}
A^{\mathsf{T}}Ax=A^{\mathsf{T}}b,
\end{equation}
which is equivalent to
\begin{equation}\label{eqI4}
\min_{x\in \mathbb{R}^n}\|b-Ax\|_2.
\end{equation}
CGLS will converge to the minimum-norm solution $x=A^\dag b$, provided $x_0\in  \rm{\mathcal{R}}(\it{A}^{\mathsf{T}})$ (See, e.g. \cite{bjorck1996numerical}, p.~291). However, the convergence of these methods deteriorates for ill-conditioned problems and they require reorthogonalization \cite{Hayami10}  to improve the convergence. Here, we say $(\ref{eq3})$ is ill-conditioned if the condition number $\kappa_2(A)=\|A\|_2\|A^\dag\|_2\gg 1$. Alternatively, the LSMR \cite{lsmr} is mathematically equivalent to applying MINRES \cite{paige1975} to $(\ref{eq4})$.\par
Hayami et al. \cite{Hayami10} proposed preconditioning the $m\times n$ rectangular matrix $A$ of the least squares problem by an
$n\times m$ rectangular matrix $B$ from the right and the left, and using the generalized minimal residual (GMRES) method \cite{saad1986} for solving the preconditioned least squares problems (AB-GMRES and BA-GMRES methods, respectively).
For ill-conditioned problems, AB-GMRES and BA-GMRES were shown to be more robust compared to the preconditioned CGNE and CGLS, respectively. Note here that BA-GMRES works with Krylov subspaces in $n$-dimensional space, whereas AB-GMRES works with Krylov subspaces in $m$-dimensional space. Since $m<n$ in the underdetermined case, AB-GMRES works in a smaller dimensional space than BA-GMRES and should be more computationally efficient compared to BA-GMRES for each iteration. Moreover, AB-GMRES has the advantage that the weight of the norm in  $(\ref{eq3})$ does not change for arbitrary $B$. Thus, we mainly focus on using AB-GMRES to solve the underdetermined least squares problem $(\ref{eq3})$. Morikuni \cite{Morikuni13} showed that AB-GMRES may fail to converge to a least
squares solution in finite-precision arithmetic for inconsistent problems. We will review this phenomenon.
The GMRES applied to inconsistent problems was also studied in other papers\cite{bw,cr,REICHEL05,Morikuni15,MorikuniRozloznik2018SIMAX}. See e.g., \cite{meza1992deflated,bw,MorikuniRozloznik2018SIMAX} for methods for solving nearly singular systems.
\par\par
In this paper, we first analyze the deterioration of convergence of AB-GMRES. To overcome the deterioration, we use the normal equations of the upper triangular matrix arising in AB-GMRES to change the inconsistent subproblem to a consistent one. In finite precision arithmetic, forming the normal equations for the subproblem will not square its condition number as would be predicted by theory. In the ill-conditioned case, the tiny singular values are shifted upwards due to rounding errors.
 Then, applying the standard Cholesky decomposition to the normal equations will result in a well-conditioned lower triangular matrix, which will ensure that the forward and backward substitutions work stably, and overcome the problem. Our approach
using the normal equations can be considered as a case where rounding errors are
beneficial \cite{nh}. We analyze why the proposed method works.
Numerical experiments on least squares problems with ill-conditioned rectangular coefficient matrices (Maragal$\_$3T to 7T \cite{florida}) show that the proposed method converges to a more accurate numerical solution than the original AB-GMRES. We also show by numerical experiments that the method is effective for applying GMRES to inconsistent range-symmetric systems with singular or severely ill-conditioned square coefficient matrices. \par
The rest of the paper is organized as follows. In Section 2, we briefly review AB-GMRES. In Section~3, we demonstrate the deterioration of convergence of AB-GMRES applied to underdetermined inconsistent least squares problems. In Section 4, we propose and present the stabilized GMRES method which is based on normal equations and has a regularization effect for ill-conditioned problems. We also explain why the method works by performing a rounding error analysis of the method. In Section 5, numerical experiment results for applying AB-GMRES to inconsistent underdetermined systems, and for applying GMRES to inconsistent systems with severely ill-conditioned and singular range-symmetric square coefficient matrices
are presented. In Section 6, we conclude the paper.\par
All the experiments in this paper were done using MATLAB R2017b in double precision, unless specified otherwise (where we extended the arithmetic precision using the Multiprecision Computing Toolbox for MATLAB \cite{mptfm}), and the computer uesd was Alienware 15 CAAAW15404JP with CPU Inter(R) Core(TM) i7-7820HK (2.90GHz).

\section{Deterioration of convergence of AB-GMRES for inconsistent problems}\label{section3AB}
In this section, we review previous work. First, we introduce the right-preconditioned GMRES (AB-GMRES). Then, we demonstrate the deterioration of convergence of AB-GMRES for inconsistent problems. Finally, we cite a related theorem to analyze the deterioration.
\subsection{AB-GMRES}
The AB-GMRES method of Hayami et al. \cite{Hayami10} applies the GMRES method \cite{saad1986} to
\begin{equation}\label{eqI5}
\min_{u\in \mathbb{R}_m}\|b-ABu\|_2,\quad x=Bu,
\end{equation}
where $B\in \mathbb{R}^{n\times m}$ is a preconditioning matrix.\par
Note the equivalence between the least squares problem (\ref{eqI4}) and the preconditioned least squares problem (\ref{eqI5}).
\begin{thm}
\rm (Theorem 3.1 of \cite{Hayami10})
\begin{equation*}
\min_{x\in \mathbb{R}^n}\|b-Ax\|_2=\min_{u\in \mathbb{R}^m}\|b-ABu\|_2
\end{equation*}
holds for all $b\in \mathbb{R}^m$ if and only if $\mathcal{R}(A)=R(AB)$.
\end{thm}

\begin{lem}\rm (Lemma 3.3 of \cite{Hayami10})\par
$\mathcal{R}(A^{\mathsf{T}})=\mathcal{R}(B)\Longrightarrow\mathcal{R}(A)=\mathcal{R}(AB)$.
\end{lem}

\begin{thm}
\rm(Theorem 3.6 of \cite{Hayami10})\it \par
If $\mathcal{R}(A^{\mathsf{T}})=\mathcal{R}(B)$, then $\mathcal{R}(AB)=\mathcal{R}(B^{\mathsf{T}}A^{\mathsf{T}}) \Longleftrightarrow \mathcal{R}(A)=\mathcal{R}(B^{\mathsf{T}})$.
\end{thm}
The convergence conditions of AB-GMRES are given as follows.
\begin{thm}\label{th3.7h}\rm
(Theorem 3.7 of \cite{Hayami10})\it \par
If $\mathcal{R}(A^{\mathsf{T}})=\mathcal{R}(B)$, then AB-GMRES determines a least squares solution of $\min_{x\in \mathbb{R}^n}\|b-Ax\|_2$ for all $b\in \mathbb{R}^m$ and for all $x_0\in \mathbb{R}^n$ if and only if $\mathcal{R}(A)=\mathcal{R}(B^{\mathsf{T}})$. Here, $x_0=Bu_0$ is the initial approximate solution of (\ref{eqI5}) when applying AB-GMRES.
\end{thm}
\par
Let $r=b-Ax=b-ABu$. Note
\begin{equation}
\|r\|_2^2=\|r|_{\rm{\mathcal{R}}(\it A)}\|_2^2+\|r|_{\rm{\mathcal{R}}(\it A)^\bot}\|_2^2=\|r|_{\rm{\mathcal{R}}(\it A)}\|_2^2+\|b|_{\rm{\mathcal{R}}(\it A)^\bot}\|_2^2.
\end{equation}
Here, $S^\bot$ denotes the orthogonal complement of a subspace $S$, and $r|_{\rm{\mathcal{R}}(\it A)}$ is the $ \rm{\mathcal{R}}(\it{A})$ component of a vector $r$. $r|_{\rm{\mathcal{R}}(\it A)^\bot}$ is the $\mathcal{R}(\it A)^\bot$ (inconsistent) componet of the residual vector $r$. Thus, AB-GMRES minimizes $\|r\|_2^2$, and hence $\|r|_{\rm{\mathcal{R}}(\it A)}\|_2^2$.\par
The $k$th iterate $x_k$ of AB-GMRES is given by
\begin{equation}
x_k=x_0+Bu_k,
\end{equation}
where $u_k\in  \mathcal{K}_k(AB, r_0)=\rm span$$\{r_0, ABr_0, \dots, (AB)^{k-1}r_0\}$, so that $x_k=x_0+z_k$, \\ where  $ z_k\in\mathcal{K}_k(BA, Br_0)=\rm span$$\{Br_0, (BA)Br_0, \dots, (BA)^{k-1}Br_0\}$.
Hence, if $x_0\in \mathcal{R}(B)$, $x_k\in \mathcal{R}(B)$.
\par
If $\mathcal{R}(B)=\mathcal{R}(A^\mathsf{T})$, then $x_k\in\mathcal{R}(A^\mathsf{T})=\mathcal{N}(A)^\bot$. Further, if $\mathcal{R}(B^\mathsf{T})=\mathcal{R}(A)$, then AB-GMRES determines a least squares solution $x_k$, i.e., $r_k|_{\rm{\mathcal{R}}(\it A)}=0$, where $r_k=b-Ax_k$,
and that solution $x_k$ is the minimum Euclidean norm solution.
\par
The algorithm is given in Algorithm \ref{AL1} \cite{Hayami10}. Here, $H_{k+1,k}=(h_{ij})\in \mathbb{R}^{(k+1)\times k}$ and $e_1=(1,0,\dots,0)^{\mathsf{T}}\in \mathbb{R}^{k+1}.$ Algorithm \ref{AB-GMRES method} is said to break down when $h_{k+1, k}=0$. See Appendix B of \cite{Morikuni15}.\par
\begin{algorithm}
\caption{AB-GMRES}
\label{AB-GMRES method}
\begin{algorithmic}[1]
\State Choose $x_0\in \mathbb{R}^{n}$, $r_0=b-Ax_0$, $v_1=r_0/\|r_0\|_2$
\For{$k=1,2,\dots$}
\State $w_k=ABv_k$
\For{$j=1,2,\dots,k$}
\State $h_{j, k}=w_k^{{\mathsf{T}}}v_j$, \quad  $w_k=w_k-h_{j, k}v_j$
\EndFor
\State $h_{k+1,k}=\|w_k\|_2$,  \quad $v_{k+1}=w_k/h_{k+1,k}$
\State Compute $y_k\in \mathbb{R}^k$ which minimizes $\|r_k\|_2=\|\|r_0\|_2\, e_1-H_{k+1,k}\, y\|_2$
\State $x_k=x_0+B[v_1, v_2, \dots, v_k]y_k$, \qquad $r_k=b-Ax_k$
\If{$\|A^{\mathsf{T}}r_k\|_2 < \epsilon\|A^{\mathsf{T}}r_0\|_2$}
\State stop
\EndIf
\EndFor
\end{algorithmic}
\label{AL1}
\end{algorithm}\par
To find $y_k\in \mathbb{R}^k$ that minimizes the $k$th residual norm $\|r_k\|_2=\|\|r_0\|_2\,e_1-H_{k+1,k}\,y_k\|_2$ in Algorithm \ref{AB-GMRES method}, the standard approach computes the QR decomposition of $H_{k+1,k}$
\begin{equation}\label{EQ3}
H_{k+1,k}=Q_{k+1}R_{k+1,k},\qquad R_{k+1,k}=\left(                 
  \begin{array}{ccc}   
    R_k\\  
    0^{\mathsf{T}}\\  
  \end{array}
\right)\in \mathbb{R}^{{(k+1)}\times{k}},
\end{equation}
where $Q_{k+1}\in \mathbb{R}^{{(k+1)}\times{(k+1)}}$ is an orthogonal matrix and $R_k\in \mathbb{R}^{{k}\times{k}}$ is an upper triangular matrix. Then, backward substitution is used to solve a system with the coefficient matrix $R_k$ as follows\par
     \begin{equation}
\|r_k\|_2=\min_{y\in \mathbb{R}^k}\| Q_{k+1}^{\mathsf{T}}(\beta e_1)-R_{k+1,k}y\|_2,
\end{equation}
where
    \begin{equation}
 \beta=\|r_0\|_2,\quad
Q_{k+1}^{\mathsf{T}}\beta e_1=\left(                 
  \begin{array}{ccc}   
    t_k\\  
    \rho_{k+1}\\  
  \end{array}
\right),\quad t_k\in \mathbb{R}^k,\quad \rho_{k+1}\in \mathbb{R}.
\end{equation}
Therefore,
 \begin{equation}
     y_k=\arg_{y\in \mathbb{R}^k}\| Q_{k+1}^{\mathsf{T}}(\beta e_1)-R_{k+1,k}y\|_2=R_k^{-1}t_k,
 \end{equation}
 and the $k$th iterate is given by
 \begin{equation}
x_k=x_0+V_ky_k,\qquad V_k=[v_1, v_2, \dots, v_k]\in \mathbb{R}^{n\times k},\qquad V_k^{\mathsf{T}}\it{V_k}=I,
\end{equation}
where $I$ is the identity matrix, and $v_1, v_2, \dots, v_k$ are the basis vectors of $\mathcal{K}_k(AB, r_0)$ defined in Algorithm \ref{AL1}.

\par
 From now on, we use AB-GMRES to solve $(\ref{eq3})$ with $B=A^{\mathsf{T}}$ and $x_0\in \mathcal{R}(A^{\mathsf{T}})$, e.g. $x_0=0$, which means $x_k=x_0+z_k$, where $z_k\in \mathcal{K}_k(A^{\mathsf{T}}A,A^{\mathsf{T}}r_0)$.
Hence, Theorem \ref{th3.7h} guarantees the convergence in exact arithmetic even in the inconsistent case.
However, in finite precision arithmetic, AB-GMRES may fail to converge to a least squares solution for inconsistent problems, as shown later.

\subsection{AB-GMRES for inconsistent problems}\label{sec23333333333}
 In this section, we perform experiments to show that the convergence of AB-GMRES deteriorates for inconsistent problems.
Experiments were done on the transpose of the matrix Maragal$\_3$  \cite{florida}, denoted by Maragal$\_3$T etc. Table \ref{tb1} gives the information on the Maragal matrices, including the density of nonzero entries, rank and condition number. Here, the rank and condition number were determined by using the MATLAB functions \texttt{spnrank} \cite{if} and \texttt{svd}, respectively.
\begin{table}
\caption{Information on the Maragal matrices.}
\begin{center}
\begin{tabular}{r|r|r|r|r|r}

 \multicolumn{1}{c|}{ matrix} & \multicolumn{1}{c|}{$m$} & \multicolumn{1}{c|}{$n$} &\multicolumn{1}{c|}{ density[$\%$]} & \multicolumn{1}{c|}{rank} & \multicolumn{1}{c}{$\kappa_2(A)$}  \\
  \hline\rule{0pt}{12pt}
  Maragal$\_$3T & 858& 1682 & 1.27 & 613& 1.10$\times 10^{3}$ \\
  Maragal$\_$4T & 1027& 1964 & 1.32 & 801&9.33$\times 10^{6}$ \\
  Maragal$\_$5T & 3296& 4654 & 0.61 & 2147& 1.19$\times 10^{5}$\\
  Maragal$\_$6T & 10144& 21251 & 0.25 & 8331 &2.91$\times 10^{6}$\\
  Maragal$\_$7T & 26525& 46845 & 0.10 & 20843&8.91$\times 10^{6}$  \\

\end{tabular}
\label{tb1}
\end{center}
\end{table}


\begin{figure}
  \centering
  \begin{tabular}{c}
  \includegraphics[width=11cm]{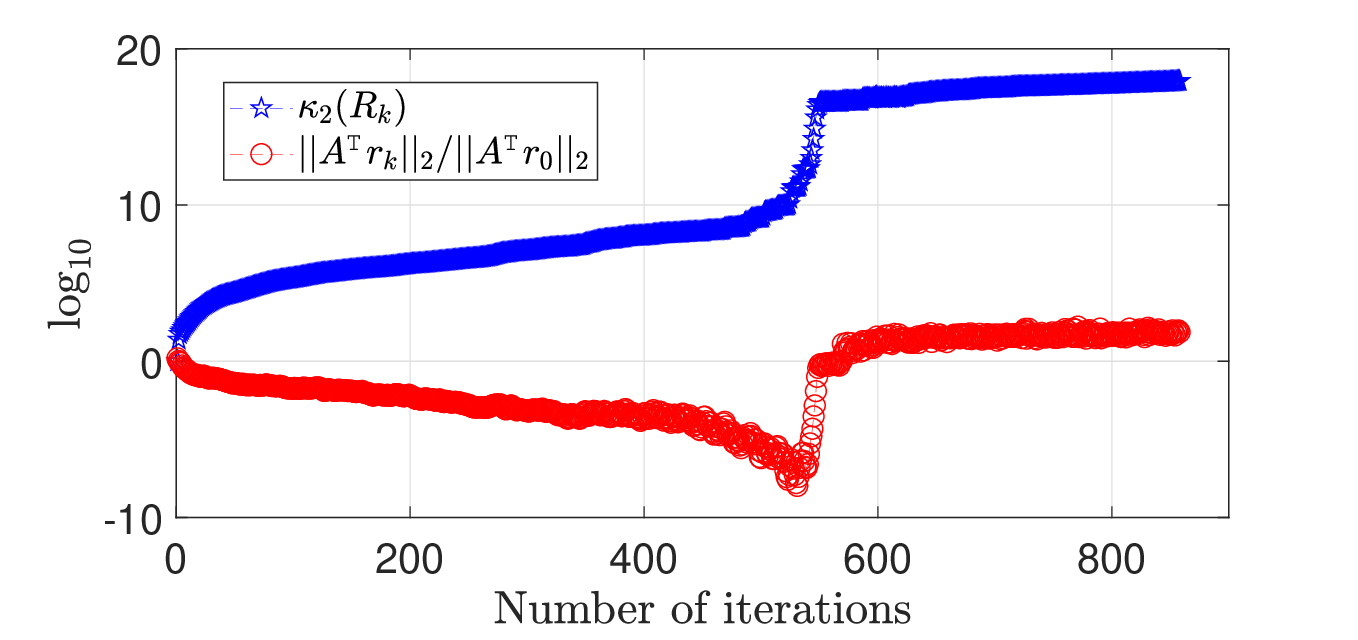}
\end{tabular}
  \caption{$\kappa_2$($R_k$) and relative residual norm versus the number of iterations for Maragal$\_$3T.}
  \label{lllw1}
 \end{figure}

Figure \ref{lllw1} shows the relative residual norm $\|A^{\mathsf{T}}r_k\|_2/\|A^{\mathsf{T}}b\|_2$ and $\kappa_2$($R_k$) versus the number of iterations for AB-GMRES with $B=A^{\mathsf{T}}$ for Maragal$\_3$T, where $r_k=b-Ax_k$, and the vector $b$ was generated by the MATLAB function \texttt{rand} which returns a vector whose entries are uniformly distributed in the interval $(0,1)$. Therefore, generically $b\notin \mathcal{R}(A)$ and the problem is inconsistent. Here, $\kappa_2$($R_k$)=$\kappa_2$($H_{k+1,k})$ holds from $(\ref{EQ3})$. The value of $\kappa_2(R_k$) was computed by the MATLAB function \texttt{cond}. The relative residual norm $\|A^{\mathsf{T}}r_k\|_2/\|A^{\mathsf{T}}b\|_2$ decreased to $10^{-8}$ until the 525th iteration, and then increased sharply. The value of \texttt{cond($R_k$)} started to increase rapidly around iterations 450--550. This observation shows that $R_k$ becomes ill-conditioned before convergence. Thus, AB-GMRES failed to converge to a least squares solution. This phenomnenon was observed by Morikuni\cite{Morikuni13}.\par
The reason why $R_k$ becomes ill-conditioned before convergence in the inconsistent case will be explained by a theorem in the next subsection.
\subsection{GMRES for inconsistent problems}\label{sec2333}
Brown and Walker \cite{bw} analyzed the break-down of GMRES.\par
Let $b|_{\rm{\mathcal{R}(\hat{\it{A}})}}$ denote the orthogonal projection of $b$ onto $\rm{\mathcal{R}(\hat{\it{A}})}$. Assume $\rm{\mathcal{N}}$$(\hat{A})=$$\rm{\mathcal{N}}$$(\hat{\it{A}}^{\mathsf{T}})$  and \\ grade$(\hat{A},b|_{\rm{\mathcal{R}(\hat{\it{A}})}})=k$. Here, grade$(\hat{A},\hat{b})$ for $\hat{A}\in \mathbb{R}^{m\times m}$, $\hat{b}\in \mathbb{R}^m$ is defined as the minimum $k$ such that $\mathcal{K}_{k+1}(\hat{A}, \hat{b})=\mathcal{K}_{k}(\hat{A}, \hat{b})$. Then,
\begin{align*}
\rm{dim}(\mathcal{K}_{\it{k}}(\hat{\it{A}}, \it b|_{\rm{\mathcal{R}(\hat{\it{A}})}}))&= \rm{dim}(\mathcal{K}_{\it k+\rm 1}(\hat{\it{A}},\it b|_{\rm{\mathcal{R}(\hat{\it{A}})}}))\\&=\rm{dim}(\hat{\it{A}}\mathcal{K}_{\it k}(\hat{\it{A}}, \it b|_{\rm{\mathcal{R}(\hat{\it{A}})}}))\\&=\rm{dim}(\hat{\it{A}}\mathcal{K}_{\it k+\rm 1}(\hat{\it{A}},\it b|_{\rm{\mathcal{R}(\hat{\it{A}})}}))\\&=k
\end{align*}
(See Appendix \ref{ap1}). Since $\rm{\mathcal{N}}$$(\hat{\it{A}})=$$\rm{\mathcal{N}}$$(\hat{\it{A}}^{\mathsf{T}})$, we obtain $\hat{A}b|_{\mathcal{R}(\hat{A})}=\hat{A}b$ and \\ 
dim($\hat{A}\mathcal{K}_{k+1}(\hat{A},b))=$dim$(\hat{A}\mathcal{K}_{k+1}(\hat{A},b|_{\rm{\mathcal{R}(\hat{\it{A}})}}))=k$. If $b\notin \mathcal{R}(\hat{A})$ and dim($\hat{A}\mathcal{K}_{k}(\hat{A},b))=k$, \\ dim($\mathcal{K}_{k+1}(\hat{A},b))=k+1$ (See Appendix \ref{ap2}).\par
Let $x_0$ be the initial solution and $r_0=b-\hat{A}x_0.$ In the inconsistent case, a least squares solution is obtained at iteration $k$, and at iteration $k+1$ breakdown occurs because of dim($\hat{A}\mathcal{K}_{k+1}(\hat{A},r_0))< $ dim($\mathcal{K}_{k+1}(\hat{A},r_0))$, i.e. rank deficiency of $\it{\min_{z\in \mathcal{K}_{k+\rm 1\it}(\hat{A},r_{\rm 0\it})}\|b-\hat{A}(x_{\rm 0\it}+z)\|_{\rm 2\it}}=\it{\min_{z\in \mathcal{K}_{k+\rm 1\it}(\hat{A},r_{\rm 0\it})}\|r_{\rm 0\it}-\hat{A}z\|_{\rm 2\it}}$\cite{bw}. This case is also called the benign breakdown\cite{REICHEL05}. \par
However, even if $\rm{\mathcal{N}}$$(\hat{A})=\rm{\mathcal{N}}$$(\hat{A}^{\mathsf{T}})$, when (\ref{eq5}) is inconsistent, the least squares problem 
\\
$\it{\min_{z\in \mathcal{K}_{k}(\hat{A},r_{\rm 0\it})}\|r_{\rm 0\it}-\hat{A}z\|_{\rm 2\it}}$ may become ill-conditioned as shown below.\par
Brown and Walker \cite{bw} introduced an effective condition number to explain why GMRES fails to converge for inconsistent least squares problems
\begin{equation}\label{eq5}
\it{\min_{x\in \mathbb{R}^m}\|b-\hat{A}x\|_{\rm 2}},
\end{equation}
where $\it{\hat{A}\in\mathbb{R}^{m\times m}}$ is singular, in the following Theorem \ref{th2}. \par

 \begin{thm}\label{th2}\rm \cite{bw}\it
\quad  Assume $\rm{\mathcal{N}}$$(\hat{A})=\rm{\mathcal{N}}$$(\hat{A}^{\mathsf{T}})$, and denote the least squares residual of \emph{(\ref{eq5})} by $r^*$, the residual at the $(k-1)$st iteration by $r_{k-1}$. If $r_{k-1}\neq r^*$, then
  \begin{equation}
  \kappa_2(A_k)\geq \frac{\|A_k\|_2}{\|\bar{A_k}\|_2}\frac{\|r_{k-1}\|_2}{\sqrt{\|r_{k-1}\|_2^2-\|r^*\|_2^2}},
  \end{equation}
where $A_k\equiv \hat{A}|_{\mathcal{K}_k(A,r_0)}$and $\bar{A_k}\equiv \hat{A}|_{\mathcal{K}_k(A,r_0)+\rm{span}\{\it{r}^*\}}$. Here, $\hat{A}|_S$ is the restriction of $\hat{A}$ to a subspace $S\subseteq \mathbb{R}^{m}$.
 \end{thm}

 Theorem \ref{th2} implies that GMRES suffers ill-conditioning for $b\notin$ $\rm{\mathcal{R}}$$(\hat{A})$ as $\|r_k\|$ approaches $\|r^*\|$.
 We can apply Theorem \ref{th2} to AB-GMRES for least-squares problems by setting $\hat{A}\equiv AA^{\mathsf{T}}$.
Theorem \ref{th2} also implies that even if we choose $B$ as $A^{\mathsf{T}}$, which satisfies the conditions in Theorem \ref{th3.7h}, AB-GMRES still may not converge numerically because of the ill-conditioning of $R_k$, losing accuracy in the solution computed in finite-precision arithmetic when $r_{k-1}$ approaches $r^*$.

\section{Deterioration of convergence of AB-GMRES applied to inconsistent underdetermined least squares problems}\label{sec3}
In this section, we illustrate, the deterioration of convergence of GMRES by numerical experiments. There are two points to note in this section. The first point is that the condition number of $R_k$ tends to become very large as the iteration proceeds for inconsistent problems, as already mentioned in section \ref{sec23333333333}. Due to $H_{k+1,k}=Q_{k+1}R_{k+1,k}$, the condition number of $H_{k+1,k}$ is the same as that of $R_k$, and will also become very large. The second point is as follows. Since $y_k=R_k^{-1}t_k$, $y_k$ is obtained by applying backward substitution to the triangular system
\begin{equation}\label{eq141l}
    R_ky_k=t_k.
\end{equation} When the triangular system becomes ill-conditioned, backward substitution becomes numerically unstable, and fails to give an accurate solution $y_k$.\par
Figure \ref{lllw1} shows that at step 550 the relative residual norm suddenly increases. To understand this increase, observe the singular values of $R_{550}$.\par
\begin{figure}
  \centering
  \begin{tabular}{c}
  \includegraphics[width=11cm]{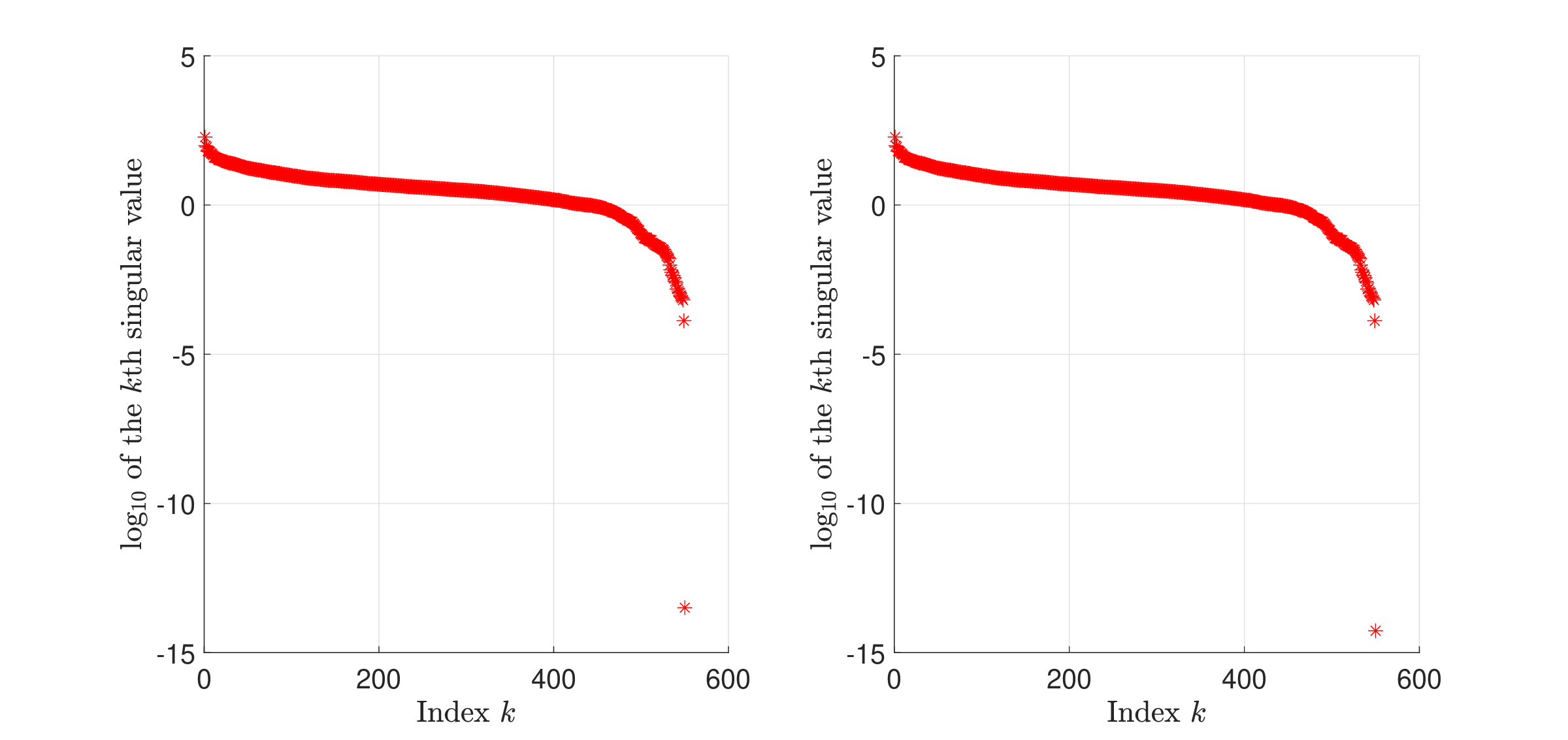}\\
  \centering{(a)Double precision arithmetic $\qquad\qquad$  (b)Quadruple precision arithmetic}
\end{tabular}
  \caption{Singular value distribution of $R_{550}$ for Maragal$\_$3T in double and quadruple precision arithmetic.}
  \label{nd}
 \end{figure}
The left of Figure \ref{nd} shows the singular values of $R_{550}$ which were computed in double precision arithmetic. The smallest singular value of $R_{550}$ is $3.21\times 10^{-14}$, which means that the triangular matrix $R_{550}$ is very ill-conditioned and nearly singular in double precision arithmetic.\par
The right of Figure \ref{nd} shows the singular values of $R_{550}$ which were computed in quadruple precision arithmetic using the Multiprecision Computing Toolbox for MATLAB \cite{mptfm}. The smallest singular value of $R_{550}$ is $5.39\times 10^{-15}$. Since quadruple precision is more accurate, from now on, we mainly show singular value distributions computed in quadruple precision.\par

 \par
 \begin{figure}
  \centering
  \begin{tabular}{c}
  \includegraphics[width=11cm]{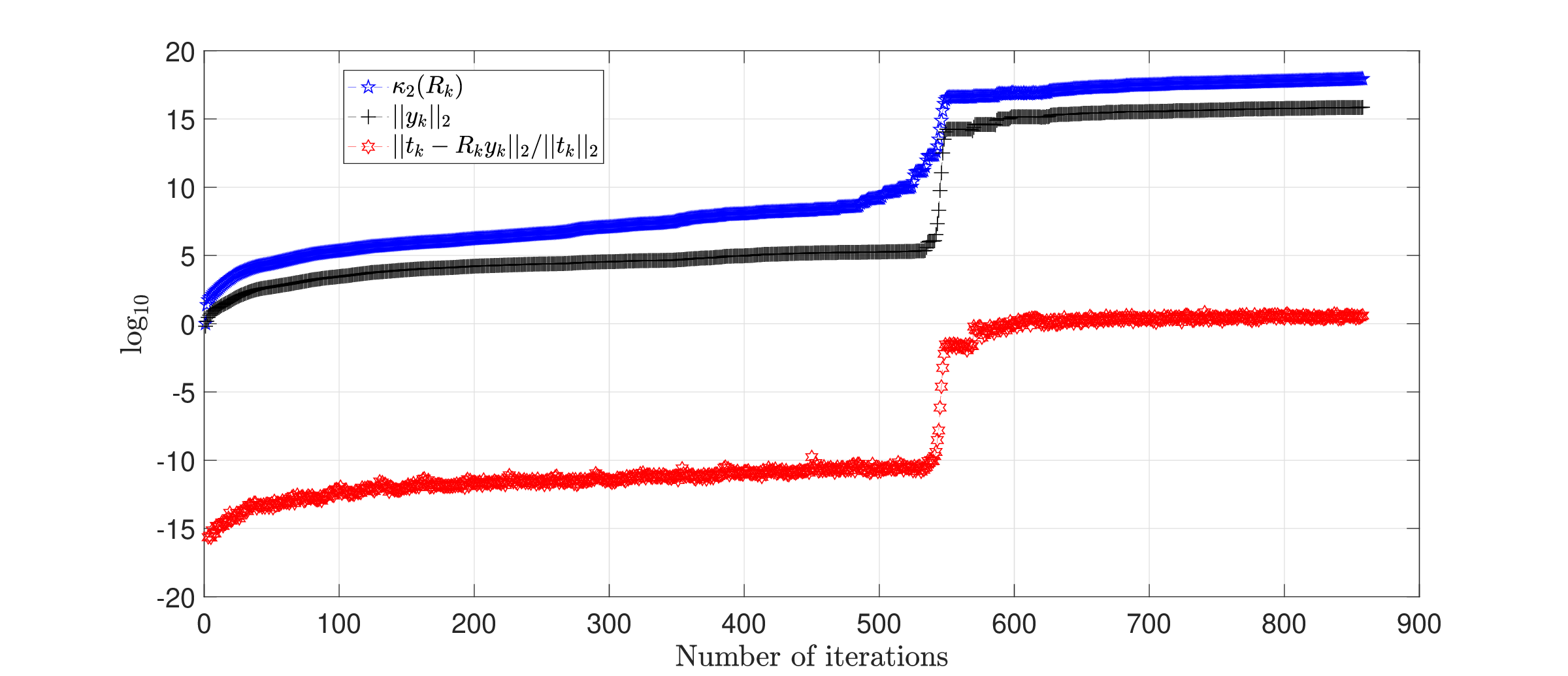}
\end{tabular}
  \caption{$\kappa_2(R_k)$, $\|y_k\|_2$, and $\|t_k-R_ky_k\|_2/\|t_k\|_2$ versus the number of iterations for Maragal$\_$3T.}
  \label{lllw5}
 \end{figure}
 Figure \ref{lllw5} shows $\kappa_2(R_k)$, $\|y_k\|_2$, and the relative residual norm $\|t_k-R_ky_k\|_2/\|t_k\|_2$ versus the number of iterations for AB-GMRES. The relative residual norm increases only gradually when the condition number of $R_k$ is less than $10^{8}$. When the condition number of $R_k$ becomes larger than $10^{10}$, the relative residual norm starts to increase sharply. This observation shows that when the condition number of $R_k$ becomes very large, the backward substitution will fail to give an accurate $y_k$. As a result, we would not get an accurate $x_k$, and the convergence of AB-GMRES would deteriorate.
\section{Stabilized GMRES method}
In this section, we first propose and present a stabilized GMRES method. Then, we explain its regularization effect comparing it with other regularization techniques.
\subsection{The stabilized GMRES}\label{sec41}
 In order to overcome the deterioration of convergence of GMRES for inconsistent systems, we propose solving the normal equations
 \begin{equation}\label{eq419}
R_k^{\mathsf{T}}R_ky_k=R_k^{\mathsf{T}}t_k
\end{equation}instead of $
R_ky_k=t_k
    $ of (\ref{eq141l}),
which we will call the stabilized GMRES.  We replace line 8 of Algorithm \ref{AL1} by Algorithm \ref{AL2}.
This makes the system consistent, and stabilizes the process, as will be shown in the following. \par
One may also consider using the normal equations of $H_{k+1,k}$. However, before breakdown, we use the standard AB-GMRES, which means we do not have to store $H_{k+1,k}$. We only store $R_k$ and update it in each iteration, which is cheaper.
\begin{figure}
  \centering
  \begin{tabular}{c}
  \includegraphics[width=11cm]{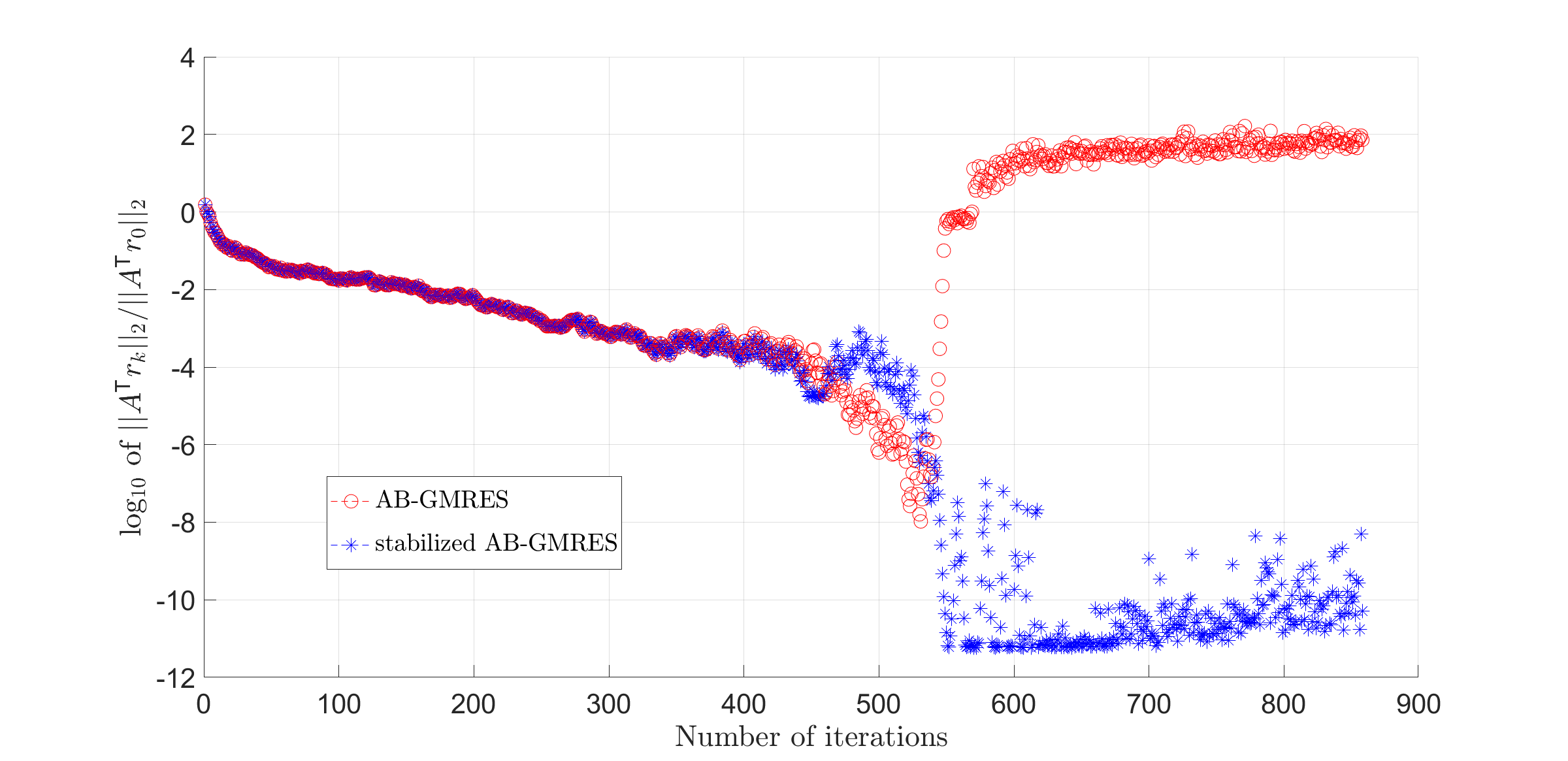}
\end{tabular}
  \caption{Comparison of the standard AB-GMRES with stabilized AB-GMRES for Maragal$\_$3T.}
  \label{lllw2}
 \end{figure}
\par
Figure \ref{lllw2} shows the relative residual norm $\|A^{\mathsf{T}}r_k\|_2/\|A^{\mathsf{T}}r_0\|_2$ versus the number of iterations for the standard AB-GMRES and stabilized AB-GMRES with $B=A^{\mathsf{T}}$ for Maragal$\_3$T. The stabilized method reaches the relative residual norm level of $10^{-11}$ which improves a lot compared to the standard method.
The method which we used for solving the normal equations (\ref{eq419}) is the standard Cholesky decomposition without pivoting.\par
\begin{algorithm}
\caption{Normal equations stabilization approach}
\label{AL2}
\begin{algorithmic}[1]
\State Compute the QR decomposition of $H_{k+1,k}=Q_{k+1}R_{k+1,k}.$
\State $R_{k+1,k}=\left(                 
  \begin{array}{ccc}   
    R_k\\  
    0^{\mathsf{T}}\\  
  \end{array}
\right)$,\qquad
 $Q_{k+1}^{\mathsf{T}}\beta e_1=\left(                 
  \begin{array}{ccc}   
    t_k\\  
    \rho_{k+1}\\  
  \end{array}
\right)$,\qquad
$\widetilde{R}_k=R_k^{\mathsf{T}}R_k,\qquad \widetilde{t}_k=R_k^{\mathsf{T}}t_k$.
\State Compute the Cholesky decomposition of $\widetilde{R}_k=LL^{\mathsf{T}}$.
\State Solve $Lz_k=\widetilde{t}_k$ by forward substitution.
\State Solve $L^{\mathsf{T}}y_k=z_k$ by backward substitution.
\end{algorithmic}
\end{algorithm}
This seems paradoxical, since forming the normal equations whose coefficient matrix $R_k^{\mathsf{T}}R_k$ would square the condition number compared to $R_k$, which would make the ill-conditioned problem even worse. Why
can the stabilized AB-GMRES give a more accurate solution? We will explain why the stabilized AB-GMRES works in the next subsection.\par

In spite of the above mentioned merits of stabilization, solving the normal equations in AB-GMRES is expensive. Actually, we only need the stabilized AB-GMRES when $R_k$ becomes ill-conditioned.
 Thus, we can speed up the process by switching AB-GMRES to stabilized AB-GMRES only when $R_k$ becomes ill-conditioned. The condition number of an incrementaly enlarging triangular matrix can be estimated by techniques in \cite{tebbens2014}.
In this paper, we adopt the switching strategy by monitoring the relative residual norm $\|A^{\mathsf{T}}r_k\|_2/\|A^{\mathsf{T}}r_0\|_2$.  Let ATR($k$)=$\|A^{\mathsf{T}}r_k\|_2/\|A^{\mathsf{T}}r_0\|_2$ for the $k$th iteration. When ATR($v$)/ $\min_{k=1, 2, \dots, v-1}$ATR($k)>10$, we judge that a jump in relative residual norm has occured, and we switch AB-GMRES to stabilized AB-GMRES at the $v$th iteration. \par

\subsection{Why the stabilized GMRES method works}
Consider solving $R_ky_k=t_k, R_k\in \mathbb{R}^{k\times k}, t_k\in \mathbb{R}^{k}$ by solving the normal equations (\ref{eq419}), which, in theory, squares the condition number and makes the problem
become harder to solve numerically. However, in 	finite precision arithmetic, the condition number of the normal equations is not neccessarily squared. We will continue to illustrate the phenomenon by using the example in Section \ref{sec3}.\par
We used the MATLAB function \texttt{svd} in quadruple precision arithmetic \cite{mptfm} to calculate the singular values. The smallest singular value of $R_{550}$ is $ 5.39\times 10^{-15}$, so its square is $2.91\times 10^{-29}$.\par
\begin{figure}
  \centering
  \begin{tabular}{c}
  \includegraphics[width=11cm]{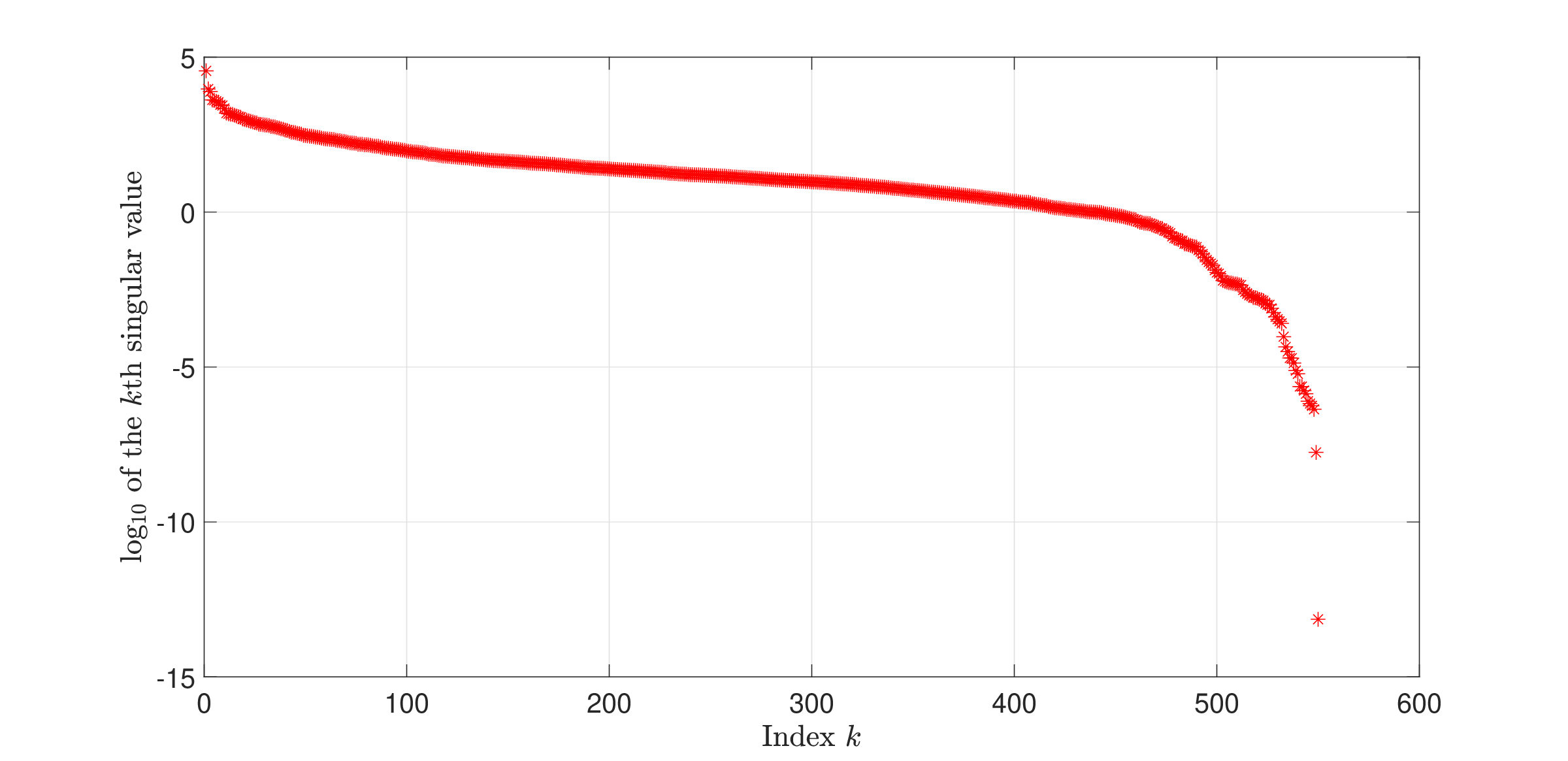}
\end{tabular}
  \caption{Singular values $\sigma_k($fl$_d(R_{550}^{\mathsf{T}}R_{550}))$, $k=1, 2, \dots, 550$ in quadruple precision arithmetic.}
  \label{lll7}
 \end{figure}
Let fl$(\cdot)$ denote the evaluation of an expression in floating point arithmetic and  fl$_d(\cdot)$ and fl$_q(\cdot)$ denote the result in double precision arithmetic and quadruple precision arithimetic, respectively.
  Figure \ref{lll7} shows that, numerically, the smallest singular value of fl$_d(R_{550}^{\mathsf{T}}R_{550})$ is $ 7.21\times 10^{-14}$, which is much larger than $2.91\times 10^{-29}$. Further, the Cholesky factor $L$ of fl$_d(R_{550}^{\mathsf{T}}R_{550})~=~LL^{\mathsf{T}}$ computed in double precision precision arithmetic has the smallest singular value $3.50\times 10^{-7}$, which is also larger than $\sqrt{2.91\times 10^{-29}}=5.39\times 10^{-15}$. Thus, the triangular systems $Lz_k=\widetilde{t}_k$ and $L^{\mathsf{T}}y_k=z_k$ are better-conditioned than $R_ky_k=t_k$, which will ensure the stability of the forward and  backward substitutions and succeeds in obtaining a much more accurate solution
 with stability compared to the standard approach as shown in Figure \ref{lllw2}. \par
  \begin{figure}
  \centering
  \begin{tabular}{c}
  \includegraphics[width=11cm]{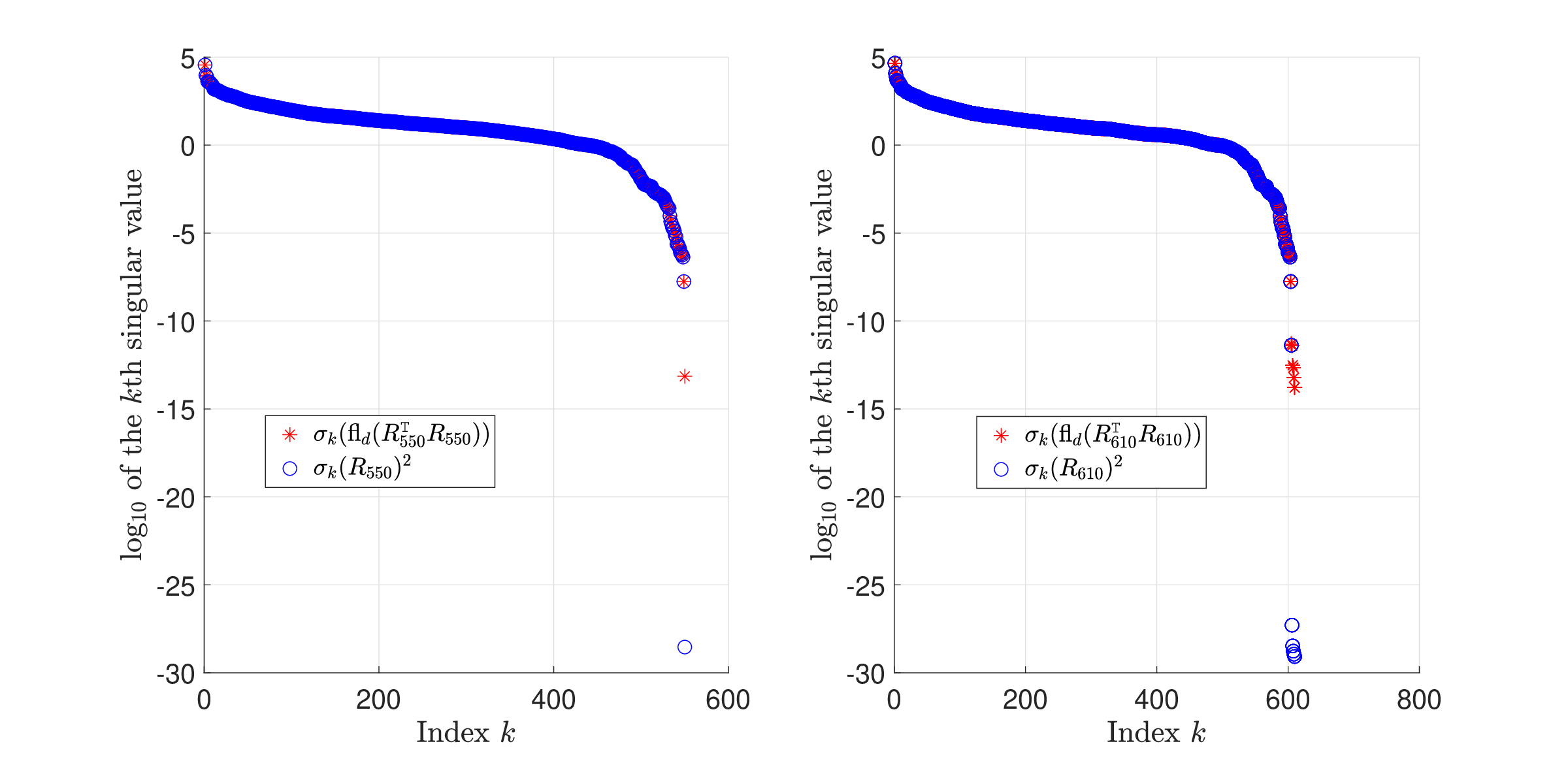}
\end{tabular}
  \caption{Singular values  $\sigma_k($fl$_d(R_{550}^{\mathsf{T}}R_{550}))$, $\sigma_k(R_{550})^2$, $\sigma_k$(fl$_d(R_{610}^{\mathsf{T}}R_{610}))$, and $\sigma_k(R_{610})^2$ in quadruple precision arithmetic.}
  \label{lll8}
 \end{figure}
 \par
The left of Figure \ref{lll8} compares the singular values $\sigma_k($fl$_d(R_{550}^{\mathsf{T}}R_{550}))$ and $\sigma_k(R_{550})^2, k=1, 2, \dots, 550$. The first to the 549th singular values of\break fl$_d(R_{550}^{\mathsf{T}}R_{550})$ and the corresponding $\sigma(R_{550})^2$ are almost the same, while the last one is different. What will happen when $R_k$ contains a cluster of small singular values?
  \par
  The upper triangular matrix $R_{610}$ contains a cluster of small singular values. The right of Figure \ref{lll8} compares the singular values  $\sigma_k($fl
$_d(R_{610}^{\mathsf{T}}R_{610}))$ and $\sigma_k(R_{610})^2$.  The larger singular values are the same as the `exact' values, while the smaller singular values become larger than the `exact' ones.
  \par
 Experiment results show that finite precision arithmetic has the effect of shifting the tiny singular values upwards and reduce the condition number of $R^{\mathsf{T}}R$. Besides the fact that the possibly inconsistent system (\ref{eq141l}) is replaced by the consistent system (\ref{eq419}), that is the reason why the normal equations (\ref{eq419}) help to make the problem easier to solve.\par
 Next, we computed $R_{550}^{\mathsf{T}}R_{550}$ in quadruple precision arithmetic and observed that the smallest singular values of $R_{550}^{\mathsf{T}}R_{550}$ coincided with the squared singular values $\sigma_k(R_{550})^2$ (blue circle symbol) in the left of Figure \ref{lll8}, unlike in double precision computation.  Since the maximum of the elements of\break $|$fl$_q(R_{550}^{\mathsf{T}}R_{550})$ $-$ fl$_d(R_{550}^{\mathsf{T}}R_{550})$ $|$ is approximately $8.16\times 10^{-12}$, double precision arithmetic contains error of the order of $10^{-12}$. Thus, double precision arithmetic has an effect of regularizing the matrix $R_{550}^{\mathsf{T}}R_{550}$, since double precision matrix multiplication is not accurate enough to keep all the information.\par
\subsection{Quadruple precision}
\begin{figure}
  \centering
  \begin{tabular}{c}
  \includegraphics[width=11cm]{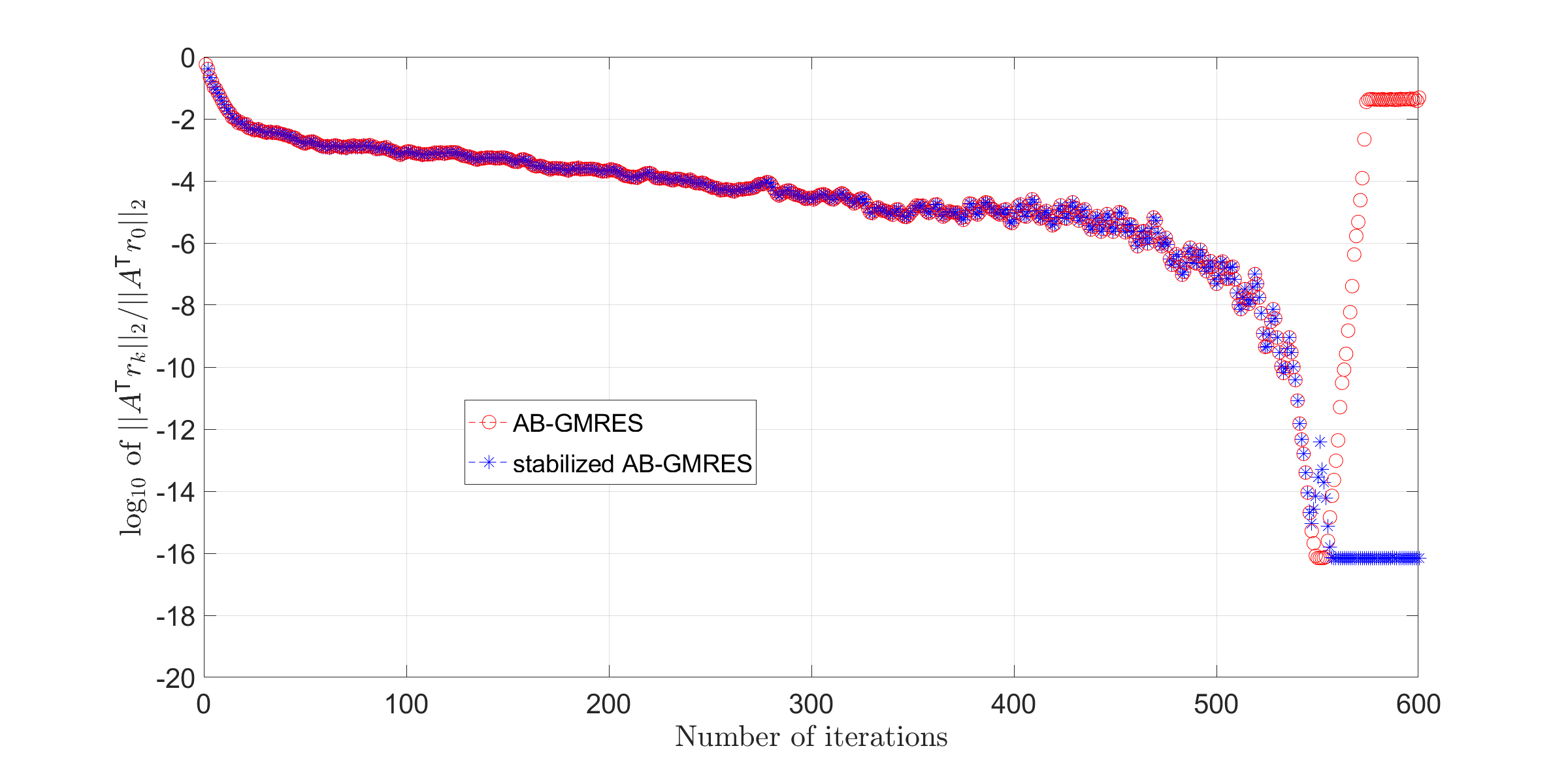}
\end{tabular}
  \caption{Effect of the stabilized method in quadruple precision arithmetic for Maragal$\_$3T.}
  \label{un1}
 \end{figure}\par
 In order to see the effect of the machine precision on the convergence of AB-GMRES, we compared the stabilized AB-GMRES with the standard AB-GMRES in quadruple precision arithmetic for the problem Maragal$\_$3T in Figure \ref{un1} in terms of the relative residual norm $\|A^{\mathsf{T}}r_k\|_2/\|A^{\mathsf{T}}b\|_2$ versus the number of iterations. For both methods, the relative residual norm reached a lower level of order $10^{-16}$ compared to $10^{-12}$ and $10^{-8}$, respectively, for double precision arithmetic in Figure \ref{lllw2}. The curves of the relative residual norm became smoother compared to double precision. As seen in Figure \ref{un1}, the relative residual norm of the standard AB-GMRES jumped to $10^{-1}$ after reaching $10^{-16}$, whereas the relative residual norm of the stabilized GMRES stayed around $10^{-16}$.

\subsection{Rounding error analysis of the stabilized GMRES method}
In order to understand the stability of the proposed method, we perform rounding error analysis. Let $u$ be the unit roundoff \cite{nh}, which is about $1.11\times 10^{-16}$ for the IEEE 754 binary 64 (double) in our experiments. The analysis shows that if the condition number of $R\in\mathbb{R}^{n\times n}$ is ${\displaystyle\frac{1}{o(n\sqrt{u})}}$, then the condition number of $R^\mathsf{T}R$ can be reduced from ${\displaystyle\frac{1}{o(n^2u)}}$ in exact arithmetic to $ O\left({\displaystyle\frac{1}{n^2u}}\right)$ in finite precision arithmetic.



Let fl$(x)$ denote the floating point number corresponding to $x$. Let $A=(a_{ij})\in\mathbb{R}^{m\times n}, B\in\mathbb{R}^{n\times p}$, and $|A|=(|a_{ij}|)$.
\par
Then, we have (cf. \cite{nh})
\begin{equation}
    |\rm f\!l\it(AB)-AB|\leq\gamma_{n}|A||B|,\quad\rm where\quad \it \gamma_n:=\frac{nu}{{\rm 1\it}-nu}.
\end{equation}
\par
Hence,
\begin{equation}
    |\rm f\!l\it(R^{\mathsf{T}}R)-R^{\mathsf{T}}R|\leq\gamma_{n}|R^{\mathsf{T}}||R|.
\end{equation}\par
Let
\begin{equation}
    E=(\varepsilon_{ij}):=\rm f\!l\it(R^{\mathsf{T}}R)-R^{\mathsf{T}}R.
\end{equation}
Then, we have
\begin{equation}
    |E|=(|\varepsilon_{ij}|)\leq\gamma_{n}|R^{\mathsf{T}}||R|.
\end{equation}

Let $R=U\Sigma V^{\mathsf{T}}$ be the singular value decomposition (SVD) of $R\in \mathbb{R}^{n\times n}$, where $U=[u_1, u_2, \dots, u_n]$,  $V=[v_1, v_2, \dots, v_n]\in\mathbb{R}^{n\times n}$ are orthogonal matrices, and $\Sigma=$ diag$(\sigma_1, \sigma_2, \dots, \sigma_n)$, $\sigma_1\geq \sigma_2\geq \cdots\geq \sigma_n\geq 0$.\par
Then, note the following.
\begin{lem}
$(|R^{\mathsf{T}}||R|)_{ij}\leq \|R\|_2^2=\sigma_1^2$.
\end{lem}
\proof\quad Let $R=[r_1 , r_2, \dots, r_n]$. Then,
\begin{align*}
	(|R^{\mathsf{T}}||R|)_{ij} & = | r_i|^\mathsf{T} |r_j| =(|r_i|,| r_j|) \\
	& \le \|  r_i \|_2 \| r_j\|_2 \le {\displaystyle \max_{1 \le i \le n} \| r_i\|_2^2} \\
	& = \max_{i=1,\dots,n} \|R e_i\|_2^2 \leq \max_{\| x\|_2=1} \|R x\|_2^2 \\
	& \le \|R\|_2^2={\sigma_1}^2.
\end{align*}
Here, $e_i$ is the $i$th column of the identity matrix. $\hfill\qed$\par

Hence, we have
\begin{equation}
    |\varepsilon_{ij}|\leq\gamma_{n}(|R^{\mathsf{T}}||R|)_{ij}\leq\gamma_{n}{\sigma_1}^2.
\end{equation}
Let  $A, B, C$ be $n\times n$ Hermitian matrices, and $A=B+C$.
Denote the eigenvalues of $A$, $B$ and $C$ by $\lambda_1(A)\geq \lambda_2(A)\geq \cdots\geq \lambda_n(A)$, $\lambda_1(B)\geq \lambda_2(B)\geq \cdots\geq \lambda_n(B)$, and $\lambda_1(C)\geq \lambda_2(C)\geq \cdots\geq \lambda_n(C)$, respectively. Then, the following Weyl's inequality (See e.g. \cite{IriM,horn2012matrix})
\begin{equation}
\lambda_k(B)+\lambda_n(C)\leq\lambda_k(A)\leq\lambda_k(B)+\lambda_1(C), \qquad  k=1, 2 \dots, n,
\end{equation}
holds.
\par
Hence, we have
\begin{equation}
|\lambda_k(A)-\lambda_k(B)|\leq\|C\|_2=\|A-B\|_2, \qquad  k=1, 2, \dots, n.
\end{equation}
Letting $A=\rm f\!l\it(R^{\mathsf{T}}R), B=R^{\mathsf{T}}R, C=E$, we have
\begin{equation}
    |\lambda_k(\rm f\!l\it(R^{\mathsf{T}}R))-\lambda_k(R^{\mathsf{T}}R)|\leq\|E\|_{\rm 2}, \qquad  k=\rm 1, 2, \dots, \it n.
\end{equation}\par
Let $\rm f\!l\it(R^{\mathsf{T}}R)=\widetilde{V}\widetilde{\Sigma}^{\rm 2}\widetilde{V}^{\mathsf{T}}$ be the SVD of $\rm f\!l\it(R^{\mathsf{T}}R)$, where $\widetilde{\Sigma}=$ diag$(\widetilde{\sigma}_1, \widetilde{\sigma}_2, \dots, \widetilde{\sigma}_n)$, $\widetilde{\sigma}_1 \geq 
\widetilde{\sigma}_2 \geq \cdots\geq \widetilde{\sigma}_n\geq 0$, which gives
\begin{equation}
    |\widetilde{\sigma}_k^2-\sigma_k^2|\leq\|E\|_2, \qquad  k=\rm 1, 2, \dots, \it n.
\end{equation}
Note
\begin{align}\label{eq241to25}
    |\widetilde{\sigma}_k^2-\sigma_k^2|&\leq\|E\|_2\leq\|E\|_F\notag\\
    &=\sqrt{\sum_{i,j=1}^n|\varepsilon_{ij}|^2}\leq\sqrt{\sum_{i,j=1}^n(\gamma_{n}{\sigma_1}^2)^2}\notag\\
    &=\sqrt{n^2(\gamma_{n}{\sigma_1}^2)^2}=n\gamma_n\sigma_1^2,\qquad  k=\rm 1, 2, \dots, \it n.
\end{align}
Hence,
\begin{equation}
    |\widetilde{\sigma}_k^2-\sigma_k^2|\leq n\gamma_n\sigma_1^2,\qquad  k=\rm 1, 2, \dots, \it n,
\end{equation}
i.e
\begin{equation}
  \sigma_k^2-n\gamma_n\sigma_1^2 \leq \widetilde{\sigma}_k^2\leq\sigma_k^2+n\gamma_n\sigma_1^2, \qquad  k=\rm 1, 2, \dots, \it n,
\end{equation}
or
\begin{equation}\label{eqsigma27}
    \widetilde{\sigma}_k^2=\sigma_k^2+t_kn\gamma_n\sigma_1^2,\quad -1\leq t_k \leq 1,\qquad  k=\rm 1, 2, \dots, \it n.
\end{equation}

Recall $\gamma_n={\displaystyle\frac{nu}{1-nu}}$, $u\approx 1.11\times10^{-16}$. If $nu\ll 1\Longleftrightarrow n\ll {\displaystyle\frac{1}{u}}$ $(\approx 9.01\times 10^{15}$ for double precision arithmetic). Then,
\begin{equation}
    {\displaystyle\frac{1}{1-nu}}\approx 1+nu                    \Longrightarrow\gamma_n={\displaystyle\frac{nu}{1-nu}}\approx nu(1+nu)\approx nu.
\end{equation}
Hence,\\
\begin{equation}
    \widetilde{\sigma}_k^2\approx\sigma_k^2+t_kn^2u\sigma_1^2,\quad -1\leq t_k \leq 1,\qquad  k=\rm 1, 2, \dots, \it n.
\end{equation}
Then,
\begin{equation}\label{leq32}
    \widetilde{\sigma}_1^2\approx\sigma_1^2(1+t_1n^2u),\qquad -1\leq t_1\leq 1.
\end{equation}\par
We define the following Landau's symbols:\par
\begin{equation}
    f(x)=O(g(x)) \text{ as } x\rightarrow \text{a denotes that }\displaystyle{\frac{f(x)}{g(x)}} \text{ is bounded as } x\rightarrow a,
\end{equation}
and
\begin{equation}
    f(x)=o(g(x))\text{ as }x\rightarrow a\text{ denotes that }\displaystyle{\lim_{x\rightarrow a}\frac{f(x)}{g(x)}=0}.
\end{equation}
In the following, for instance, $o(n\sqrt{u})$ is defined by letting $x=n\sqrt{u}, a=0.$\par
Assume $n\sqrt{u}\ll 1$ $(\Longleftrightarrow n\ll {\displaystyle\frac{1}{\sqrt{u}}}\approx 9.49\times 10^7)$. Then, since, $\displaystyle{\frac{1}{\sqrt{u}}\ll\frac{1}{u}}$, we have $n\ll \displaystyle{\frac{1}{\sqrt{u}}\ll\frac{1}{u}}$. Thus, $nu\ll 1$. \par Assume
\begin{equation*}\label{con34}
    \kappa=\kappa(R)=\frac{\sigma_1}{\sigma_n}=\frac{1}{o(n\sqrt{u})}
\end{equation*}
\begin{equation}
 \Longleftrightarrow   \kappa^{\rm 2\it}=\frac{\rm 1\it }{o(n^2u)}\quad \rm \Longleftrightarrow \quad  \frac{1}{\kappa^{\rm 2}}\it={o(n^{\rm 2}\it u)}
\end{equation}
holds. Then,
\begin{align}
    \widetilde{\sigma}_n^2&\approx\sigma_n^2+t_nn^2u\sigma_1^2=\sigma_1^2\left(\frac{\sigma_n^2}{\sigma_1^2}+t_nn^2u\right)\notag\\
    &=\sigma_1^2\left(\frac{1}{\kappa^2}+t_nn^2u\right)\approx\sigma_1^2\ [o(n^2u)+t_nn^2u].\label{leq37}
\end{align}\par
Assume $|t_n|>o(1)$. (Note that if we assume $t_n$ is randomly distributed in the interval $[-1, 1]$, then, generically, $o(1)<|t_n|$ holds.) Then, since $-1\leq t_n\leq 1$, $o(1)<|t_n|=O(1)$ holds.
Hence, $\widetilde{\sigma}_n^2\approx \sigma_1^2t_nn^2u$. Since $\widetilde{\sigma}_n^2\geq 0$, we have $t_n>0$ and ${\displaystyle\frac{1}{t_n}=O(1)}$.\par
Hence, from (\ref{leq32}), (\ref{leq37}) and $\displaystyle{\frac{1}{t_n}}=O(1)$, we have
\begin{equation}
   \widetilde{\kappa}^2=\kappa(\rm f\!l\it(R^{\mathsf{T}}R))=\frac{\widetilde{\sigma}_{\rm 1}^{\rm 2\it}}{\widetilde{\sigma}_n^{\rm 2}}\approx\frac{\rm 1}{\it t_n n^{\rm 2}\it u}\approx O\left(\frac{\rm1}{\it n^{\rm2}\it u}\right).
\end{equation}

In summary, we have the following theorem.
\begin{thm}\label{thm7l}
Let $u$ be the unit roundoff, and $R\in \mathbb{R}^{n\times n}$. If $n\sqrt{u}\ll 1$ and $\kappa(R)={\displaystyle\frac{\sigma_1(R)}{\sigma_n(R)}}={\displaystyle\frac{1}{o(n\sqrt{u})}}$, then, generically, $\sigma_n(\rm f\!l$$(R^\mathsf{T}R))\approx \sigma_1(R)^2\,t_n \,n^2u$, where $o(1)<t_n=O(1)$, and $\kappa(\rm f\!l\it(R^{\mathsf{T}}R))= O\left({\displaystyle\frac{\rm 1}{\it n^{\rm 2\it}u}}\right)$ hold.
\end{thm}
\noindent
{\bf Remark 1} \quad
For IEEE double $u\approx 1.11\times 10^{-16}$, $n\sqrt{u}\ll 1\Longleftrightarrow n\ll 9.49\times 10^{7}$. \par
\vspace{3mm}
\noindent
{\bf Remark 2} \quad
$\kappa(R^{\mathsf{T}}R)={\displaystyle\frac{1}{o(n^2u)}}$. \par
\vspace{3mm}
\noindent
Then, numerical experiments suggest that if $LL^{\mathsf{T}}$ is the Cholesky decomposition of $\rm f\!l\it(R^{\mathsf{T}}R)$ computed in finite precision, then $\kappa(L)=O\left({\displaystyle\frac{\rm 1}{\it n\sqrt{u}}}\right)$, even when $\kappa(R)={\displaystyle\frac{1}{o(n\sqrt{u})}}$.\par

Thus, forming the normal equations and applying Cholesky decomposition can lead to a more stable computation for extremely ill-conditioned systems of equations, and hence explains why the stabilized GMRES method works without choosing the value of a regularization parameter such as in TSVD or Tikhonov regularization, which will be mentioned in § \ref{sec511} and \ref{sec512}, respectively.

Let us compare estimates with numerical results for the Maragal$\_$3T matrix in Figure \ref{lll8}. \par
For $R_{550}$, $n=550$, $\sigma_{1}(R_{550})\approx 1.90\times 10^{2}$, $\sigma_{550}(R_{550})\approx 5.39\times^{-15}$. Hence,
\begin{equation*}
    \kappa(R_{550}) = \frac{\sigma_{1}(R_{550})}{\sigma_{550}(R_{550})}\approx  3.53         \times10^{16}=\frac{1}{o(n\sqrt{u})}\gg \frac{1}{n\sqrt{u}}
\approx 1.73\times 10^{5}.
\end{equation*}
Thus, $\widetilde{\sigma}_{n}^{2} \approx \sigma_1^2\,t_n\,n^2u\approx 1.21\times 10^{-6}$, where $o(1)<t_n=O(1)$, and\break $\kappa(\rm f\!l\it(R^{\mathsf{T}}R))\approx O\left({\displaystyle\frac{\rm1}{\it n^{\rm2}\it u}} \right) \approx\rm 2.98\times 10^{10}$, whereas in Figure \ref{lll8}, $\widetilde{\sigma}_n^{\rm 2}\approx 7.21\times 10^{-14}$, and $\kappa(\rm f\!l\it(R^{\mathsf{T}}R))\approx\rm 5.00\times 10^{17}$.\par
 For $R_{610}$, $n=610$, $\sigma_{1}(R_{610})\approx 2.13\times 10^{2}$, $\sigma_{610}(R_{610})\approx 2.91\times 10^{-15}$. Hence,
\begin{equation*}
    \kappa(R_{610}) = \frac{\sigma_{1}(R_{610})}{\sigma_{610}(R_{610})}\approx 7.32\times 10^{16}=\frac{1}{o(n\sqrt{u})}\gg \frac{1}{n\sqrt{u}}
\approx 1.56\times 10^{5}.
\end{equation*}
Thus, $ \widetilde{\sigma}_{n}^{2} \approx \sigma_1^2\,t_n\,n^2u\approx 1.87\times 10^{-6}$, where $o(1)<t_n=O(1)$, and\break $\kappa(\rm f\!l\it(R^{\mathsf{T}}R))\approx O\left({\displaystyle\frac{\rm1}{\it n^{\rm2}\it u}}\right)\approx\rm 2.42\times 10^{10}$, whereas in Figure \ref{lll8}, $\widetilde{\sigma}_n^{\rm 2}\approx 1.62\times 10^{-14}$, and $\kappa(\rm f\!l\it(R^{\mathsf{T}}R))\approx\rm 2.77\times 10^{18}$.\par
We summarize the results in Table \ref{tbrrrrrr}. We think there are two reasons for the overestimation of $\Tilde{\sigma}_n^2$. One comes from the inequality $\|E\|_2\leq\|E\|_F$ in (\ref{eq241to25}). The other is that $t_n>o(1)$, but $t_n$ may be considerably samller than 1 in (\ref{eqsigma27}).
\begin{table}
\caption{Comparison of estimates and numerical experiments for Maragal$\_$3T }
\begin{center}
\setlength{\tabcolsep}{1mm}
\begin{tabular}{c|r|r|r|r}

 \multirow{2}{*}{}&  \multicolumn{2}{c|}{$R_{550}$ $(n=550)$} & \multicolumn{2}{c}{$R_{610}$ $(n=610)$}\\ \cline{2-5}\rule{0pt}{12pt}
 &\multicolumn{1}{c|}{$\widetilde{\sigma}_n^2$} &\multicolumn{1}{c|}{$\widetilde{\sigma}_1^2$/$\widetilde{\sigma}_n^2$ } & \multicolumn{1}{c|}{$\widetilde{\sigma}_n^2$} & \multicolumn{1}{c}{$\widetilde{\sigma}_1^2$/$\widetilde{\sigma}_n^2$}\\
\hline\rule{0pt}{12pt}
Estimates    &1.21$\times 10^{-6}$                & 2.98$\times 10^{10}$                  &  1.87$\times 10^{-6}$             &  2.42$\times 10^{10}$               \\ \hline\rule{0pt}{12pt}
Numerical experiment (Figure \ref{lll8})   &$7.21\times10^{-14} $                &$5.00\times 10^{17} $                &  $1.62\times 10^{-14} $            & $2.77\times 10^{18}$     \\
\hline\rule{0pt}{12pt}
\multirow{2}{*}{} &\multicolumn{1}{c|}{$\sigma_n^2$} &\multicolumn{1}{c|}{$\sigma_1^2$/$\sigma_n^2$ } & \multicolumn{1}{c|}{$\sigma_n^2$} & \multicolumn{1}{c}{$\sigma_1^2$/$\sigma_n^2$}\\ \cline{2-5}\rule{0pt}{12pt}
 &\multicolumn{1}{c|}{ $2.91\times 10^{-29} $} &\multicolumn{1}{c|}{ $1.25\times 10^{33} $ } & \multicolumn{1}{c|}{ $8.47\times10^{-30} $} & \multicolumn{1}{c}{ $5.36\times 10^{33} $}
                 \end{tabular}
\end{center}
\label{tbrrrrrr}
\end{table}\par
We remark that \cite{yamamoto2015roundoff} analyzes the stability of the CholeskyQR2 algorithm using similar techniques. However, they assume $\kappa(R)\leq O\left({\displaystyle\frac{1}{\sqrt{u}}}\right)$, whereas we assume $\kappa(R)= {\displaystyle\frac{1}{o(n\sqrt{u})}}.$

\subsection{Two advantages of forming the normal equations}
When $R$ is singular, $R^{-1}$ does not exist, and
\begin{equation}\label{eq38l}
    Ry=t
\end{equation}
does not have a solution when $t\notin \mathcal{R}(R)$.\par
If we reformulate (\ref{eq38l}) as a least squares problem
\begin{equation}\label{eq39l}
    \min_{y}\|t-Ry\|_2,
\end{equation}
then (\ref{eq39l}) has a solution even when $t\notin \mathcal{R}(R)$. For instance, the minimum-norm solution of (\ref{eq39l}) is given by $y=R^{\dag}t$, where $R^{\dag}$ is the pseudo-inverse of $R$.\par
Note that (\ref{eq39l}) is equivalent to the normal equations
\begin{equation}\label{eq40l}
    R^{\mathsf{T}}Ry=R^{\mathsf{T}}t.
\end{equation}
(\ref{eq40l}) is consistent, i.e. $R^{\mathsf{T}}t\in \mathcal{R}(R^{\mathsf{T}})=\mathcal{R}(R^{\mathsf{T}}R)$, and has a solution.\par
Now consider the case when $R$ is nearly singular (severely ill-conditioned), that is $\kappa(R)=\displaystyle{\frac{1}{o(n\sqrt{u})}}$. Then, solving (\ref{eq38l}) by backward substitution fails to give an accurate solution as shown in Figure \ref{lllw5}.\par
On the other hand, we may expect that we may obtain a numerical solution of the least squares problem (\ref{eq39l}), for instance by approximating $y=R^{\dag}t$ \cite{shz}.\par
In fact, since (\ref{eq39l}) is equivalent to the normal equations, as we have seen in Theorem \ref{thm7l}, $\kappa(\rm f\!l\it(R^{\mathsf{T}}R))= O\left({\displaystyle\frac{\rm 1}{\it n^{\rm 2\it}u}}\right)$ holds while $\kappa(R^{\mathsf{T}}R)={\displaystyle\frac{1}{o(n^2{u})}}$, which gives a numerical advantage.\par
Thus, we may say that forming the normal equations (\ref{eq40l}) has two advantages over the system of equations (\ref{eq38l}). One is that, it makes the system consistent and guarantees the existence of a solution, which opens the possibility of a numerical solution by some kind of approximation. The other advantage is that normal equations become numerically better conditioned than in exact arithmetic.
 \subsection{When the stabilized GMRES method works}{\label{sec43}}
The stabilized GMRES does not always stabilize the solution of the upper triangular system. A counter example is when $R_k$ is a L\"{a}uchli matrix \cite{higham}, implying that $R^{\mathsf{T}}_kR_k$ computed in finite precision becomes singular. Indeed, when GMRES is applied to a linear system
with an EP (equal preojecton) matrix $A_3$, that is $\rm{\mathcal{N}}$$(A_3)$$=$$\rm{\mathcal{N}}$$(A_3^{\mathsf{T}})$ such as
\begin{equation}\label{exep}
  A_3x=\left(
            \begin{array}{ccccc}
               \frac{\sqrt{2}}{2} &\qquad & \frac{\sqrt{2}}{2}- \frac{\sqrt{6u}}{6} &\qquad&- \frac{\sqrt{6u}}{6}\\
                \frac{\sqrt{2}}{2}&\qquad & \frac{\sqrt{2}}{2}+ \frac{\sqrt{6u}}{6}&\qquad& \frac{\sqrt{6u}}{6} \\
               0 &\qquad& \frac{\sqrt{6u}}{3} &\qquad&\frac{\sqrt{6u}}{3}\\
             \end{array}
           \right)x= \left(
             \begin{array}{c}
1\\
0\\
0\\
   \end{array}
           \right),
\end{equation}
where  $A_3$ has the null space $\rm{\mathcal{N}}$$(A_3)=\rm{span}\{(1,-1,1)^\mathsf{T}\}$, and $u$ is the unit roundoff, the resulting $R_k$ is a L\"{a}uchli matrix.\par
 Apply GMRES with $x_0=0$ to (\ref{exep}). Let $R_k\in \mathbb{R}^{k\times k}$ be the upper triangular matrix obtained at the $k$th iteration of GMRES. In the second iteration, after applying the Givens rotation to $H_{3, 2}$, we obtain the following:
\begin{equation}\label{excep1}
R_2= \left(
             \begin{array}{cc}
               1 & 1 \\
               0 & \sqrt{u} \\
             \end{array}
           \right),\qquad R_2^{\mathsf{T}}R_2= \left(
             \begin{array}{cc}
               1 & 1 \\
               1 & 1+u \\
             \end{array}
           \right)\simeq \left(
             \begin{array}{cc}
               1 & 1 \\
               1 & 1\\ \end{array}
           \right).
\end{equation}
Thus, there is a risk that the stabilized GMRES will give a numerically singular matrix $R_2^{\mathsf{T}}R_2$ in finite precision arithmetic for nonsingular $R_2$. We will analyze this phenomenon.\par

Note that the following theorem holds from Theorem 8.10 of \cite{nh}, where $|b|=(|b_1|, |b_2|, \dots, |b_n|)^{\mathsf{T}}$ for $b=(b_1, b_2, \dots,b_n)^{\mathsf{T}}\in \mathbb{R}^{n}$.
\begin{thm}\label{th333}
Let $T=(t_{ij})\in \mathbb{R}^{n\times n}$ be a triangular matrix and $b\in \mathbb{R}^{n}.$ Then, the computed solution $\hat{x}$ obtained from substitution applied to $Tx=b$ satisfies
\begin{equation}
\hat{x} = x+O(n^2u)M(T)^{-1}|b|.
\end{equation}
Here, $M(T)=(m_{ij})$ is the comparison matrix such that
\end{thm}
\begin{equation}
   m_{ij}=\left\{
             \begin{array}{cc}
               |t_{ij}|, & i=j, \\
               -|t_{ij}|, & i\neq j. \\
             \end{array}
           \right.
\end{equation}\par
\smallskip
Further, we define the following. Let
 \begin{equation}\mathbb{O}(x)=\left(  \begin{array}{c}
                O(x)\\
                O(x)\\
                \vdots\\
              O(x) \\
             \end{array}\right)\in \mathbb{R}^n ,\quad \mathcal{O}(x)= [\mathbb{O}(x), \mathbb{O}(x), \cdots, \mathbb{O}(x)]\in \mathbb{R}^{n\times n}.\end{equation}
We assume that the basic arithmetic operations op $=$ $+, -, *, /$ satisfy\break fl$(x$ op $y)= (x $ op $y)(1+O(u))$ as in
\cite{nh}.\par Let $x, y\in \mathbb{R}^n$, $A\in \mathbb{R}^{n\times n}$. Then,\par
\begin{center}
    fl$(x^\mathsf{T}y)=x^\mathsf{T}y+O(nu)|x|^{\mathsf{T}}|y|=x^{\mathsf{T}}y+O(nu)$,\par
fl$(Ax)=Ax+\mathbb{O}(nu)|A||x|=Ax+\mathbb{O}(nu)$.
\end{center}\par
Let $C\in \mathbb{R}^{n\times n}$ and $\|C\|_2=O(1)$. We say $C\in \mathbb{R}^{n\times n}$ is numerically nonsingular if the statement
\begin{equation}\label{EQ14}
\rm f\!l\it (Cx)=\mathbb{O}(u)\quad\Rightarrow\quad x=\mathbb{O}(u)
\end{equation}
holds. Note that this definition of numerical nonsingularity agrees with that of numerical rank \cite{bjorck1996numerical} due to the following.\par
Let the SVD of $C=U\Sigma V^{\mathsf{T}}$, where $U, V$ are orthogonal matrices and $\Sigma=\rm diag\it (\sigma_{\rm 1}, \sigma_{\rm 2}, \dots, \sigma_n).$ We assume $\|C\|_2=\sigma_1=O(1).$ If the numerical rank of $C$ is $r<n$, there is a singular value $\sigma_i=O(u),$ $r+1\leq i\leq n.$ Then, $Cx=U\Sigma V^{\mathsf{T}}x=\mathbb{O}(u)$ admits $x'=V^{\mathsf{T}}x=(x_1', x_2', \dots, x_n')^{\mathsf{T}}$ such that $x_i'=O(1)$, and hence $x=\mathbb{O}(1).$ Thus, $C$ is numericaly singular. Then, the following theorem holds.
\begin{thm}\label{th55}
Let $R_k=(r_{ij})\in \mathbb{R}^{k\times k}$ be an upper-triangular matrix and
\begin{equation}
   R_{k+1}=\left(
             \begin{array}{cc}
               R_k & d \\
              0^{\mathsf{T}} & r_{k+1,k+1} \\
             \end{array}
           \right)\in \mathbb{R}^{(k+1)\times(k+1)}.\end{equation}
 Assume that $R_k$ is nonsingular and numerically nonsingular, $R_k=\mathcal{O}(1), \break R_k^{-1}=\mathcal{O}(1), M(R_k)^{-1}=\mathcal{O}(1), d=\mathbb{O}(1)$, and $O(k)=O(k^2)=O(1).$ Then, the following holds:
\begin{equation*}
\rm f\!l\it(R_{k+{\rm 1\it}}^{\mathsf{T}}R_{k+{\rm 1\it}})~is~numerically~nonsingular \it\ \Longleftrightarrow\ \rm f\!l \it (r^{\rm 2\it}_{k+{\rm 1\it},k+{\rm 1\it}})> \rm f\!l \it(d^{\mathsf{T}}d)O(u).
\end{equation*}
\end{thm}
\begin{proof} See Appendix \ref{ap3}. \end{proof}
Theorem \ref{th55} gives the necessary and sufficient condition so that the stabilized GMRES works at the $(k+1)$st iteration, i.e. $ R_{k+1}^{\mathsf{T}}R_{k+1}$ is numerically nonsingular.
\begin{figure}
  \centering
  \begin{tabular}{c}
  \includegraphics[width=11cm]{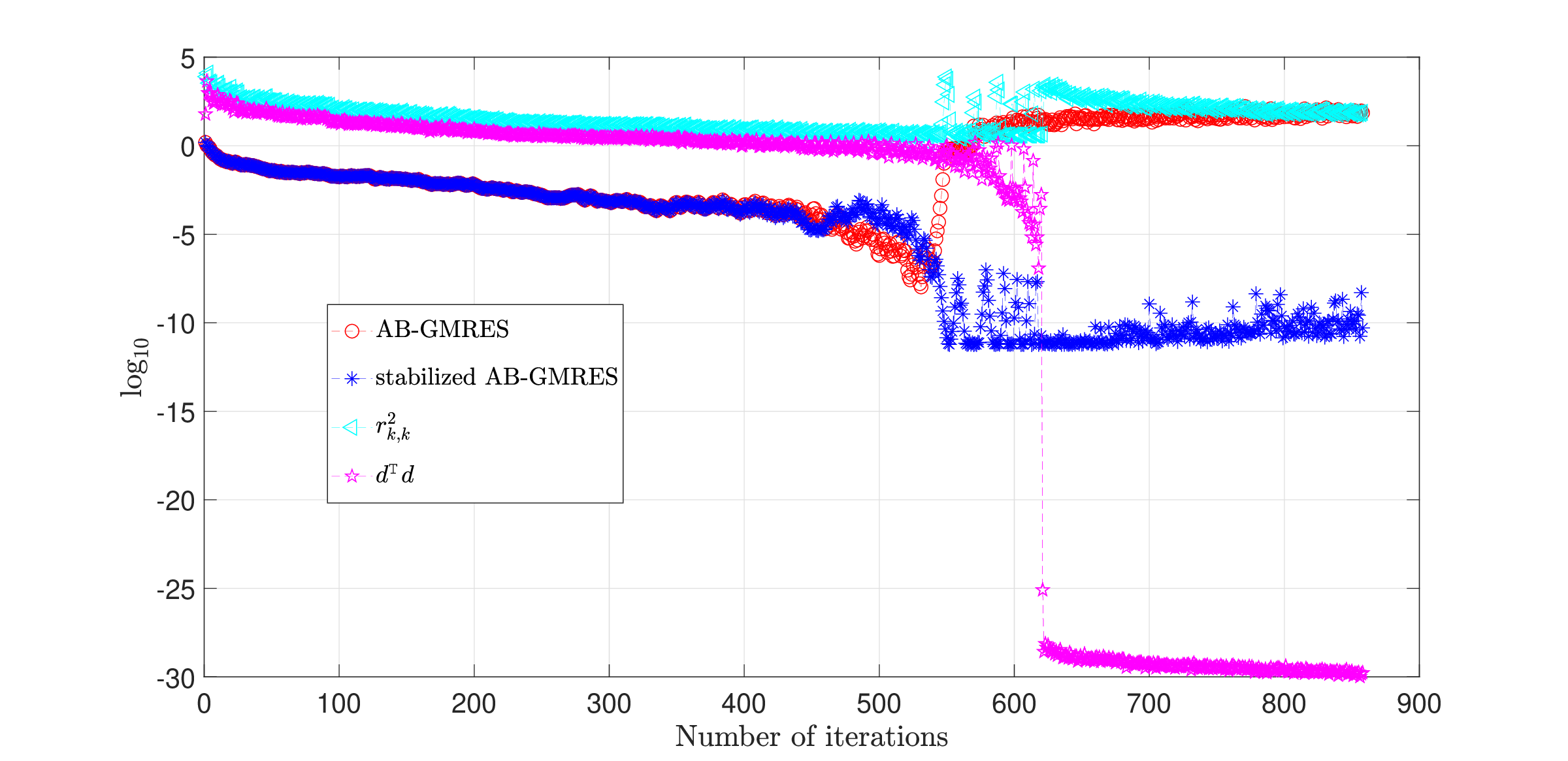}
\end{tabular}
  \caption{$r_{k,k}^2$, $d^{\mathsf{T}}d$, and $\|A^{\mathsf{T}}r_k\|_2/\|A^{\mathsf{T}}b\|_2$ for AB-GMRES and stabilized AB-GMRES for Maragal$\_$3T.}
  \label{t4dr}
 \end{figure}\par

The difficulty in solving $R_iy_i=t_i$ by backward substitution is not necessarily because the diagonals of $R_i$ are tiny. The reason is that $R_i$ has  tiny singular values. However, the exceptional example (\ref{excep1}) exists where the stabilized AB-GMRES does not work. The condition fl$(r^2_{k+1, k+1})$ $>$ fl$(d^{\mathsf{T}}d)O(u)$ in Theorem \ref{th55} excludes such exceptions.\par
Figure \ref{t4dr} shows $r^2_{k+1, k+1}$ and $d^{\mathsf{T}}d$ together with the relative residual norm $\|A^{\mathsf{T}}r_k\|_2/\|A^{\mathsf{T}}b\|_2$ of AB-GMRES and stabilized AB-GMRES for Maragal$\_$3T. The figure shows that up to 613 iterations, the conditions in Theorem \ref{th55} are satisfied, and $ R_{k+1}^{\mathsf{T}}R_{k+1}$ is numerically nonsingular, so that the stabilized AB-GMRES works.\par
In fact, for
\begin{equation*}
R_2= \left(
             \begin{array}{cc}
               1 & 1 \\
               0 & \sqrt{u} \\
             \end{array}
           \right)
\end{equation*}
  of  (\ref{excep1}), $ \sigma_1(R_2)\approx\sqrt{2}$, $\sigma_2(R_2)\approx {\displaystyle\sqrt{\frac{u}{2}}} $, so that $ \kappa(R_2)\approx {\displaystyle\frac{2}{\sqrt{u}}}\approx O\left({\displaystyle\frac{1}{n\sqrt{u}}}\right)\ll o\left({\displaystyle\frac{1}{n\sqrt{u}}}\right)$, so the condition (\ref{con34}) is not satisfied, and the stabilized GMRES is not guaranteed to work in this case.

\section{Comparisons with other methods}
We show the numerical performance of the proposed stabilized AB-GMRES method on test matrices, compared with previous methods. All programs for iterative methods were coded according to the algorithms in \cite{NEUMAN2,Hayami10,lsqr,lsmr}.  Each method was terminated at the iteration step which gives the minimum relative residual norm within $m$ iterations, where $m$ is the number of the rows of the matrix. No restarts were used for GMRES. Experiments were done for rank-deficient underdetermined matrices whose information is given in Table~1. Here, we have deleted the zero rows and columns of the test matrices beforehand. The elements of $b$ were randomly generated using the MATLAB function \texttt{rand}. Therefore, generically $b\notin \mathcal{R}(A)$ and the problem is inconsistent. Each experiment was done 10 times for the same right-hand side $b$ and the average of the CPU times are shown. The symbol - denotes that  $\|A^{\mathsf{T}}r_k\|_2/\|A^{\mathsf{T}}r_0\|_2$ did not reach $10^{-8}$ within $m$ iterations. The symbol $(*)$ denotes that we used the MATLAB function \texttt{chol} instead of Cholesky decomposition without pivoting for solving the normal equations (\ref{eq419})  to save CPU time, except for Havard500, for which Cholesky decomposition without pivoting did not converge. The symbol $(\&)$ denotes the case where even using the MATLAB function \texttt{chol} for solving equation (\ref{eq419}) failed to  converge. Then, we used the MATLAB function backslash for solving the normal equations (\ref{eq419}).

\subsection{Underdetermined inconsistent least squares problems} 

\subsubsection{Comparison with Truncated SVD method}\label{sec511}
Motivated by the stabilized AB-GMRES, we also applied the truncated singular value decomposition (TSVD) stabilization method and compared it with the stabilized AB-GMRES. The method modifies $R_k$ by truncating singular values smaller than~$\mu$.
More specifically, let $R_{k}=U\Sigma V^{\mathsf{T}}$ be the SVD of $R_k$, where the columns of $U=[u_1, u_2, \dots, u_k]$ and $V=[v_1, v_2, \dots, v_k]$ are the left and right singular vectors, respectively, and the diagonal entries of $\Sigma=$ diag$(\sigma_1, \sigma_2, \dots, \sigma_k)$ are the singular values of $R_k$ in discending order $\sigma_1\geq \sigma_2\geq \cdots\geq \sigma_k$. Then, the TSVD approximates $R_k\simeq \sum_{i=1}^j \sigma_i u_i v_i^{\mathsf{T}}$ with $j$ such that $\sigma_{j+1}\leq \mu\sigma_1\leq \sigma_j$ and $y_k= R_k^{-1}t_k\simeq \sum_{i=1}^j\frac{1}{\sigma_i}v_i u_i^{\mathsf{T}}t_i, j\leq k$.\par
\begin{figure}
  \centering
  \begin{tabular}{c}
  \includegraphics[width=11cm]{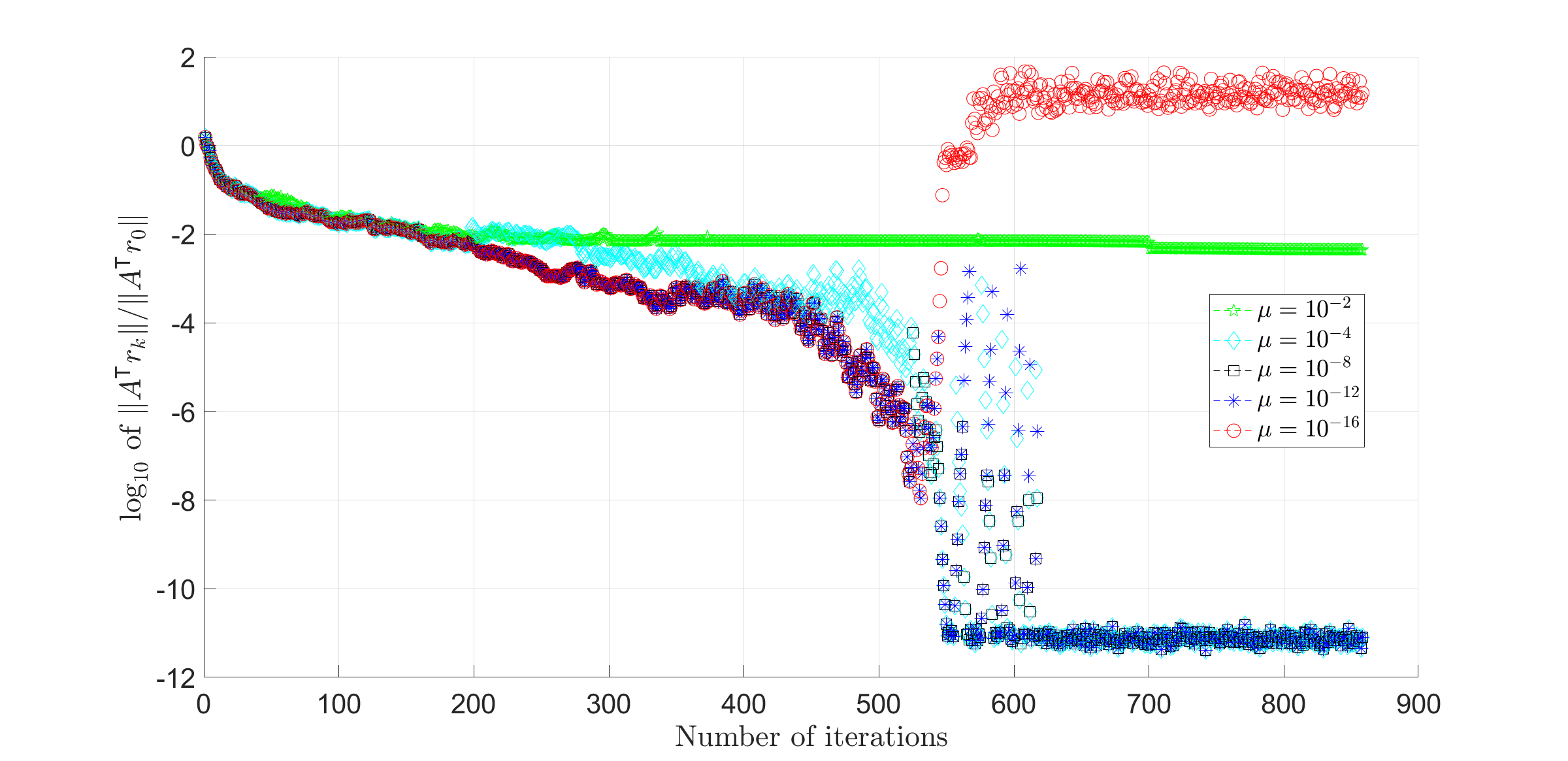}
\end{tabular}
  \caption{Relative residual norm for TSVD stabilized AB-GMRES versus number of iterations for different values of the regularization parameter $\mu$ for Maragal$\_$3T.}
  \label{lwp13}
 \end{figure}\par
\begin{figure}
  \centering
  \begin{tabular}{c}
  \includegraphics[width=11cm]{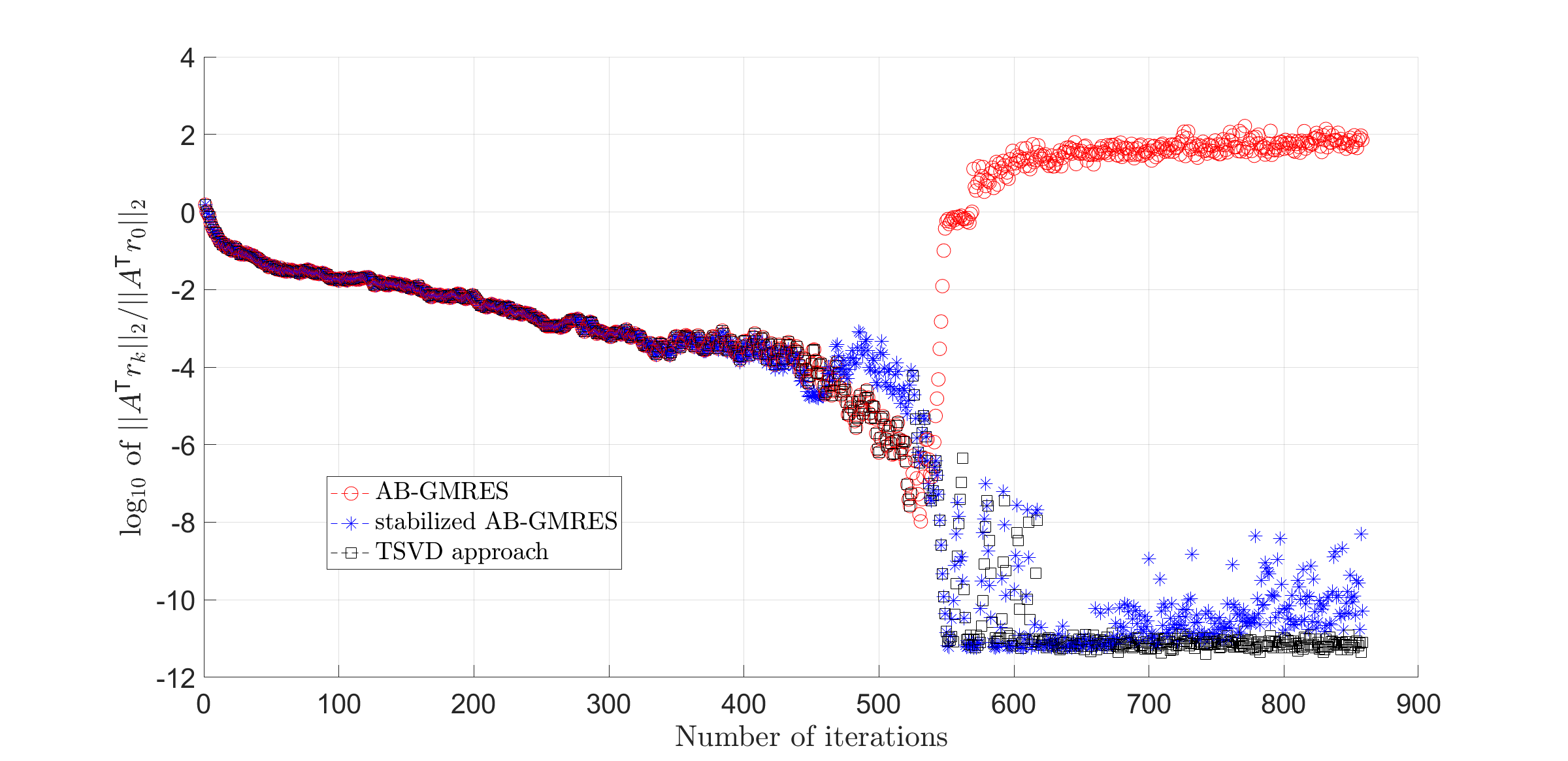}
\end{tabular}
  \caption{Comparison of the standard AB-GMRES with stabilized and TSVD stabilized AB-GMRES with $\mu=10^{-8}$ for Maragal$\_$3T.}
  \label{lwp4}
 \end{figure}
\par
For the problem Maragal$\_$3T, when $\mu= 10^{-13}, 10^{-12}, \dots, 10^{-4}$, the method converges but when $\mu$ is smaller than $10^{-13}$ or larger than $10^{-4}$, it does not converge as shown in Figure \ref{lwp13}. Numerical experiments showed that $\mu=\sqrt{u}\simeq 10^{-8}$, where $u$ is the unit roundoff (about $10^{-16}$ in double presion arithmetic), gave the best result among $\mu = 10^{-1}, 10^{-2},\dots, 10^{-16}$ in terms of the relative residual norm. The  convergence behaviour of  the TSVD stabilization method with $\mu=10^{-8}$ is similar to the stabilized AB-GMRES method as shown in Figure \ref{lwp4}, which suggests that eliminating tiny singular values of $R_k$ which are less than $10^{-8}$ is effctive for sovling problem (\ref{eq3}). However, the TSVD method requires computing the truncated singular value decomposition of $R_k$, and requires choosing the value of the threshold parameter $\mu$, whereas the stabilized AB-GMRES does not require either of them.

\subsubsection{Comparison with Tikhonov regularization method}\label{sec512}
Another approach to stabilize AB-GMRES would be to apply Tikhonov regularization. There are two methods to implement it. The first method is to solve the following square system:
\begin{equation}\label{eq66}
(R_k^{\mathsf{T}}R_k+\lambda I)y_k=R_k^{\mathsf{T}}t_k,\qquad \lambda\geq 0
\end{equation}
using the Cholesky decomposition. The second method is to solve the regularized least suqares problem
\begin{equation}\label{eqz16}
  \min_{y_k\in\mathbb{R}^k}\left|\left|  \left( \begin{array}{c}
t_k\\
0\\
   \end{array}\right)-\left(\begin{array}{c}
R_k\\
\sqrt{\lambda} I\\
   \end{array} \right)y_k\right|\right|_2
\end{equation}
using the QR decomposition. These two methods are equivalent  mathematically. However, they are not equivalent  numerically. The behavior of the first method is similar to the stabilized AB-GMRES.\par Table \ref{tb33} shows that AB-GMRES combined with the first method converges better when $\lambda=10^{-16}$ than when $\lambda=10^{-14}$ for the problem Maragal$\_$3T. This method can be used to shift upwards the small singular values, but is less acurrate compared to the stabilized AB-GMRES. 
\\
\begin{figure}
  \centering
  \begin{tabular}{c}
  \includegraphics[width=11cm]{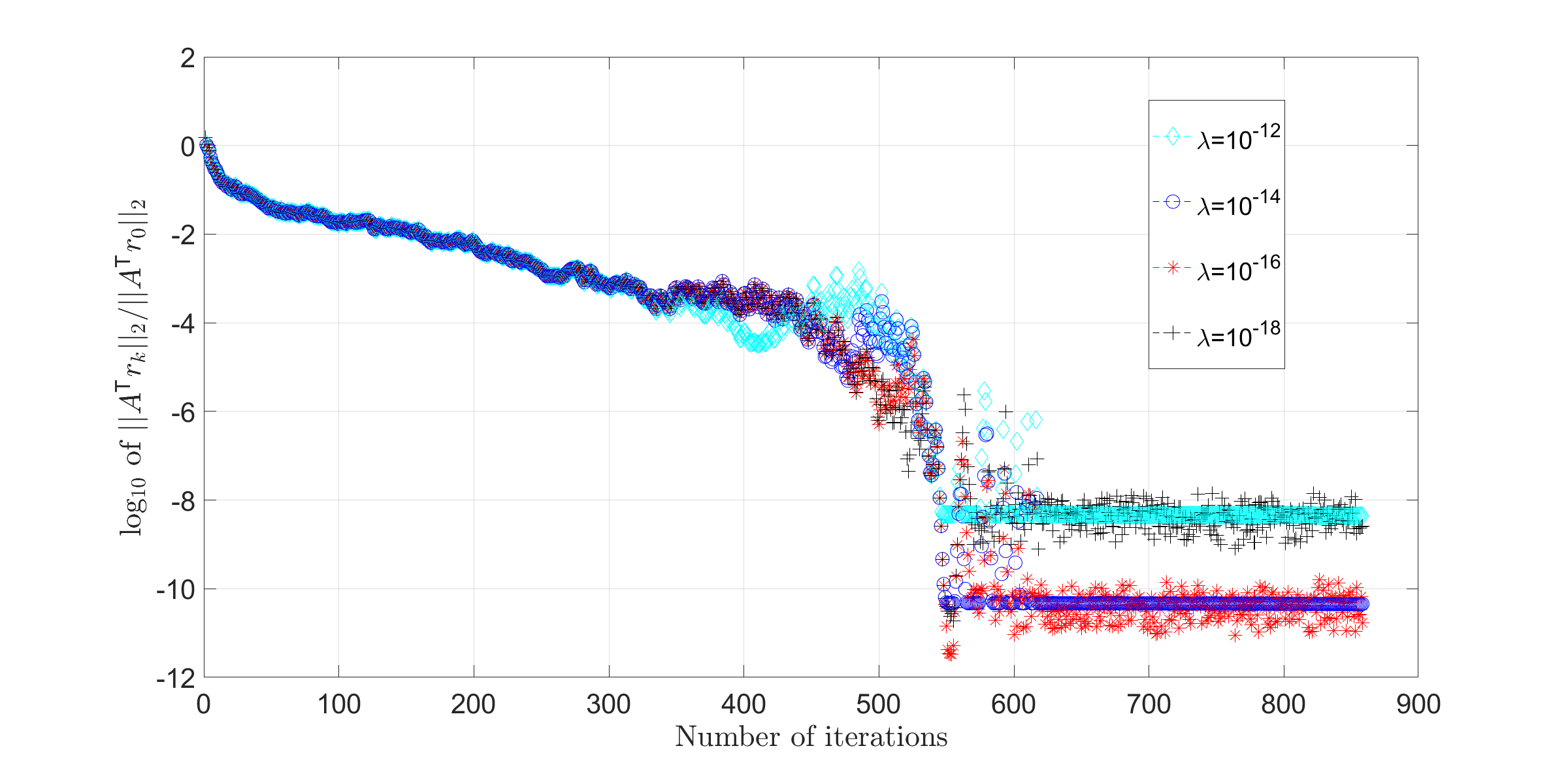}
\end{tabular}
  \caption{Relative residual norm for AB-GMRES with Tikhonov regularization using (\ref{eqz16}) versus number of iterations for different values of the regularization parameter $\lambda$ for Maragal$\_$3T.}
  \label{llwre999}
 \end{figure}\par
Table \ref{tb33} also shows that the second method is even more accurate compared with the stabilized AB-GMRES method. There is no need to form the normal equations, so that less information is lost due to rounding error. 
However, one needs to choose an appropriate value for the regularization parameter $\lambda$. Figure \ref{llwre999} shows the relative residual norm $\| A^\mathsf{T} r_k \|_2 / \| A^\mathsf{T} r_0\|_2$ for AB-GMRES  with Tikhonov regularization  using     (\ref{eqz16}) versus the number of
iterations for different values of $\lambda$ for Maragal$\_$3T. According to Figure~\ref{llwre999},    $\lambda=10^{-16}$ was optimal among $10^{-12}, 10^{-14}, 10^{-16}$, and $10^{-18}$.\par 
We here note the following.
\begin{thm}\label{thmeighth}
Let $\sigma_1\geq \sigma_2\geq\cdots\geq\sigma_k$ be the the singular values  of $R_k$. Then, the singular values of \begin{equation}R_k'=\left( \begin{array}{c}
R_k\\
\sqrt{\lambda}I\\
   \end{array}\right)\end{equation} are given by $\sqrt{\sigma_1^2+\lambda}\geq\sqrt{\sigma_2^2+\lambda}\geq\cdots\geq\sqrt{\sigma_k^2+\lambda}.$
\end{thm}
\begin{proof} See Appendix \ref{ap4}. \end{proof}
\smallskip\par
Then, let \begin{equation}\kappa\equiv\kappa_2(R_k)=\frac{\sigma_1}{\sigma_k},\quad \kappa'^2\equiv\kappa_2(R_k')^2=\frac{\sigma_1^2+\lambda}{\sigma_1^2/\kappa^2+\lambda}=1+\frac{\sigma_1^2(1-1/\kappa^2)}{\sigma_1^2/\kappa^2+\lambda}.\end{equation} Since $\kappa\geq 1,  \rm d\it \kappa'/\rm d\it \lambda \leq \rm 0$ for $\lambda\geq 0$ and $\kappa'(\lambda=0)=\kappa,  \kappa'(\lambda=+\infty)=1.$ Note also that
\begin{equation}
\lambda=\frac{\sigma_1^2[1+(\kappa'/\kappa)^2]}{\kappa'^2-1}.
\end{equation}
Therefore, for instance, if $\kappa\gg 1$ and we want $\kappa'=\sqrt{\kappa}$,
\begin{equation}
\lambda=\frac{\sigma_1^2(1+1/\kappa)}{\kappa-1}\simeq \frac{\sigma_1^2}{\kappa}.
\end{equation}
For example, if $\kappa=10^{16}$ and we want $\kappa'=10^8$, we should choose $\lambda\simeq \sigma_1^2\times 10^{-16}.$ For Maragal$\_$3T, the largest singular value~$\sigma_1$ is about 12.64, so
that we can estimate a reasonable value of $\lambda \simeq 1.60 \times 10^{-14}$. However, this estimation assumes $\kappa'=\sqrt{\kappa}$, and needs an extra cost for computing $\sigma_1$. See \cite{brezinski2009error} for other estimation
techniques for the regularization parameter.
\begin{table}
\caption{Attainable smallest relative residual norm $\|A^{\mathsf{T}}r_k\|_2/\|A^{\mathsf{T}}r_0\|_2$ for AB-GMRES with Tikhonov regularization using (\ref{eq66}) and (\ref{eqz16}), and stabilized AB-GMRES for Maragal$\_$3T.}
\begin{center}
\setlength{\tabcolsep}{0.3mm}
\begin{tabular}{c|rrrrr}

matrix & Maragal$\_$3T & Maragal$\_$4T & Maragal$\_$5T & Maragal$\_$6T & Maragal$\_$7T\\
\hline\rule{0pt}{12pt}
iter.    &552                  & 597                  &  1304              &  2440                &  1864\\
method (\ref{eq66})  $\lambda=10^{-14}$  &5.08$\times 10^{-11}$   &  5.57$\times 10^{-8}$       &  1.05$\times 10^{-5}$     &  8.26$\times 10^{-6}$        &  4.53$\times 10^{-6}$\\
\hline\rule{0pt}{12pt}
iter.    &570                  & 598                  &  1226              &  2440                & 1864\\
method (\ref{eq66})  $\lambda=10^{-16}$    &5.80$\times 10^{-12}$   &  5.59$\times 10^{-8}$      & 4.22$\times 10^{-6}$        &  8.26$\times 10^{-6}$&  4.53$\times 10^{-6}$\\

\hline\rule{0pt}{12pt}
iter.    &553                  & 547                  &  1261              &  2937                &  2475\\
method (\ref{eqz16})  $\lambda=1.6\times 10^{-14}$    &7.54$\times 10^{-11}$   &  5.59$\times 10^{-8}$      & 1.15$\times 10^{-5}$        &  9.12$\times 10^{-6}$&  2.78$\times 10^{-7}$\\
\hline\rule{0pt}{12pt}
iter.    &551                  & 547                  &  1262              &  3037                &  2475\\
method (\ref{eqz16})  $\lambda=10^{-16}$    &3.37$\times 10^{-12}$   &  5.59$\times 10^{-8}$      & 5.64$\times 10^{-7}$        &  1.91$\times 10^{-6}$&  2.78$\times 10^{-7}$\\

\hline\rule{0pt}{12pt}
iter.    &552                  &($\&$) 598                   &($*$)  1224               &($*$)  3000                 &($*$)  2475 \\
stabilized AB-GMRES    &4.86$\times 10^{-12}$   &  5.59$\times 10^{-8}$      & 2.54$\times 10^{-6}$        &  4.56$\times 10^{-6}$&  2.78$\times 10^{-7}$\\
                 \end{tabular}
\end{center}
\label{tb33}
\end{table}
\subsubsection{Comparison with the Range Restricted GMRES}
We compared the proposed stabilized AB-GMRES with the range restricted AB-GMRES (RR-AB-GMRES) \cite{NEUMAN2}, where the Krylov subspace for the RR-AB-GMRES with $B=A^{\mathsf{T}}$ is $\mathcal{K}_k(AA^{\mathsf{T}}, AA^{\mathsf{T}}r_0)$, and the standard AB-GMRES with $B=A^{\mathsf{T}}$.
\begin{table}
\caption{Comparison of the attainable smallest relative residual norm $\|A^{\mathsf{T}}r_k\|_2/\|A^{\mathsf{T}}r_0\|_2$.}
\begin{center}
\setlength{\tabcolsep}{0.5mm}{
\begin{tabular}{c|rrrrr}

matrix & Maragal$\_$3T & Maragal$\_$4T & Maragal$\_$5T & Maragal$\_$6T & Maragal$\_$7T\\
\hline\rule{0pt}{12pt}
iter.    &531                  & 465                  &  1110              &  2440                &  1864\\
standard AB-GMRES    &1.05$\times 10^{-8}$   &  2.09$\times 10^{-7}$       &  5.35$\times 10^{-6}$     &  8.26$\times 10^{-6}$        &  4.53$\times 10^{-6}$\\
\hline\rule{0pt}{12pt}
iter.    &552                  &($\&$) 598                   &($*$)  1224               &($*$)  3000                 &($*$)  2475 \\
stabilized AB-GMRES    &4.86$\times 10^{-12}$   &  5.59$\times 10^{-8}$      & 2.54$\times 10^{-6}$        &  4.56$\times 10^{-6}$&  2.78$\times 10^{-7}$\\
\hline\rule{0pt}{12pt}
iter.    &553                & 565                  &  1223              &  2374                &  2474\\
RR-AB-GMRES    &2.57$\times 10^{-11}$   &  5.59$\times 10^{-8}$      & 3.62$\times 10^{-6}$        &  1.63$\times 10^{-5}$&  2.78$\times 10^{-7}$\\
                 \end{tabular}}
\end{center}
\label{tb22}
\end{table}
\par

Table \ref{tb22} gives the number of iterations and the smallest relative residual norm for the RR-AB-GMRES, the standard and stabilized AB-GMRES for the Maragal matrices. The table shows that the stabilized AB-GMRES is more accurate than the standard AB-GMRES. Table \ref{tb22} also shows that the stabilized AB-GMRES is generally more accurate than the RR-AB-GMRES.  The stabilized AB-GMRES took more iterations to attain the same order of the smallest residual norm than the RR-AB-GMRES.\par
 \subsection{Inconsistent systems with severely ill-conditioned range-symmetric coefficient matrices}
Next, we test the stabilized AB-GMRES on least squares problems
$\min_{x\in \mathbb{R}^n} \|b-Ax\|_2$ by GMRES, where $A\in \mathbb{R}^{n\times n}$ are severely ill-conditioned range-symmetric (square) matrices
given in Table \ref{tb94}. \par
These matrices are all numerically singular.
We generated the right-hand side $b$ by the MATLAB function \texttt{rand}, so that the systems are generically inconsistent. We compared the stabilized AB-GMRES with the standard AB-GMRES and RR-AB-GMRES.
Table \ref{tb95} gives the smallest relative residual norm and the corresponding number of iterations. Table \ref{tb695} gives the CPU times in seconds required to obtain relative residual norm $\|A^{\mathsf{T}}r_k\|_2/\|A^{\mathsf{T}}r_0\|_2<10^{-8}$. The switching strategy which was introduced in Section \ref{sec41} was used for the stabilized AB-GMRES when measuing CPU times. The number of iterations when switching occurred is in brackets.
\par
\begin{table}
\caption{Information of the singular square matrices.}
\begin{center}
\setlength{\tabcolsep}{0.5mm}
\begin{tabular}{r|r|r|r|r|r}

  \multicolumn{1}{c|}{matrix} & \multicolumn{1}{c|}{size} & \multicolumn{1}{c|}{density[$\%$]} & \multicolumn{1}{c|}{rank}  & \multicolumn{1}{c|}{$\kappa_2(A)$} &\multicolumn{1}{c}{application}\\
  \hline\rule{0pt}{12pt}
  Harvard500 & 500 & 1.05 &170 &1.30$\times 10^2$ & web connectivity \\
netz4504 & 1961 & 0.13  & 1342&3.41$\times 10^1$ & 2D/3D finite element problem\\
  TS & 2142 & 0.99  & 2140 &3.52$\times 10^3$ &counter example problem\\
  grid2$\_$dual & 3136 & 0.12 & 3134 &8.58$\times 10^3$  &2D/3D finite element problem\\
uk  & 4828 & 0.06 & 4814 & 6.62$\times 10^3$&undirected graph \\
bw42  & 10000 & 0.05 & 9999 &2.03$\times 10^{3}$ & partial differential equation\cite{bw} \\
msc01050  & 1050 & 2.38 & 1049 &1.31$\times 10^{8}$ & 2D/3D structural problem \\
freeFlyingRobot$\_$7 & 3918 & 0.20 & 3881 &1.68$\times 10^{12}$ & optimal control problem \\

\end{tabular}
\end{center}
\label{tb94}
\end{table}
\begin{table}
\caption{Comparison of the attainable smallest relative residual norm $\|A^{\mathsf{T}}r_k\|_2/\|A^{\mathsf{T}}r_0\|_2$ for inconsistent square linear systems.}
\begin{center}
\setlength{\tabcolsep}{1mm}

\begin{tabular}{c|rrrrrr}

matrix & Harvard500 & netz4504 &TS & grid2$\_$dual & uk &  bw42\\
\hline\rule{0pt}{12pt}
iter.    &104                  & 144                 &  1487              &  3134                &  4620   &715\\
standard\\ AB-GMRES    &9.38$\times 10^{-9}$   &  4.51$\times 10^{-10}$       &1.56$\times 10^{-9}$     &5.98$\times 10^{-10}$&  1.35$\times 10^{-9}$        &  8.06$\times 10^{-8}$\\
\hline\rule{0pt}{12pt}
iter.    &($*$) 134                  &($\&$) 201                 &  1613              & ($*$) 3135                 &($*$)  4739 & ($\&$) 788 \\
stabilized\\ AB-GMRES    &8.46$\times 10^{-14}$   & 1.51$\times 10^{-14}$      & 2.51$\times 10^{-9}$        &  5.53$\times 10^{-10}$&6.57$\times 10^{-10}$&  1.66$\times 10^{-7}$\\
\hline\rule{0pt}{12pt}
iter.    &135                & 200                  &  1652              &  3134                &  4706&1163\\
RR-\\AB-GMRES    &7.78$\times 10^{-14}$   &  3.36$\times 10^{-14}$      & 4.56$\times 10^{-9}$        &  6.52$\times 10^{-8}$&  8.33$\times 10^{-8}$&1.56$\times 10^{-5}$\\
                 \end{tabular}
\end{center}
\label{tb95}
\end{table}
\par


\begin{table}
\caption{Comparison of the CPU time (seconds) to obtain relative residual norm $\|A^{\mathsf{T}}r_k\|_2/\|A^{\mathsf{T}}r_0\|_2<10^{-8}$ for inconsistent square linear systems.}
\begin{center}
\setlength{\tabcolsep}{0.60mm}
\begin{tabular}{c|rrrrrr}

matrix & Harvard500 & netz4504 &TS & grid2$\_$dual & uk &  bw42\\
\hline\rule{0pt}{12pt}
iter.    &104                  & 134                 &  1411              &  3134                &  4583   &-\\
standard AB-GMRES    &4.72$\times 10^{-2}$   &  1.87$\times 10^{-1}$       &2.14$\times 10$     &2.16$\times 10^{2}$&  6.93$\times 10^{2}$        &  -\\
\hline\rule{0pt}{12pt}
iter.    &104                 & 134                &  1531 (182)             &  3134                &  4679 (4199)& -\\
stabilized AB-GMRES    &4.78$\times 10^{-2}$   & 1.89$\times 10^{-1}$      & 8.19$\times 10$        &  2.21$\times 10^{2}$&1.93$\times 10^{3}$&  -\\
\hline\rule{0pt}{12pt}
iter.    &114                & 153                  &  1530              &  -                &  -&-\\
RR-AB-GMRES    &6.42$\times 10^{-2}$   &  2.62$\times 10^{-1}$      & 2.68$\times 10$        &  -&  -&-\\
                 \end{tabular}
\end{center}
\label{tb695}
\end{table}

\par
\begin{table}
\caption{Attainable smallest relative residual norm $\|A^{\mathsf{T}}r_k\|_2/\|A^{\mathsf{T}}r_0\|_2$ for range symmetric matrices.}
\begin{center}
\begin{tabular}{c|rrrrrr}

matrix & bw42 & msc01050 & freeFlyingRobot$\_$7\\
\hline\rule{0pt}{12pt}
iter.    &147                  & 560                 &  1084  \\
standard GMRES    & 8.08$\times 10^{-9}$   &  4.98$\times 10^{-8}$       &8.86$\times 10^{-8}$    \\
\hline\rule{0pt}{12pt}
iter.    &219                 & 668                &  3414     \\
stabilized GMRES    &2.11$\times 10^{-11}$   & 4.62$\times 10^{-9}$      & 3.24$\times 10^{-10}$   \\
\hline\rule{0pt}{12pt}
iter.    &220                & 564                  &  3183          \\
RR-GMRES    &3.13$\times 10^{-11}$   &  2.62$\times 10^{-6}$      & 1.40$\times 10^{-9}$  \\
                 \end{tabular}
\end{center}
\label{tbbw42}
\end{table}

\par


For Harvard500 and bw42, AB-GMRES could only converge to the level of $10^{-9}$ regarding the relative residual norm, while the stabilized AB-GMRES converged to the level of $10^{-14}$. The stabilized AB-GMRES was robust in the sense that it could continue to compute even when the upper triangular matrix $R_k$ became seriously ill-conditioned, and the relative residual norm did not increase sharply towards the end, but just stagnated at a low level, just like for consistent problems.

Thus, our stabilization method also makes AB-GMRES stable for highly ill-conditioned inconsistent systems with square coefficient matrices.\par
The coefficient matrix $A$ of bw42 is singular and satisfies $\rm{\mathcal{N}}$$(A)=$$\rm{\mathcal{N}}$$(A^{\mathsf{T}})$. The problem comes from a finite-difference discritization of a PDE with periodic boundary condition (Experiment 4.2 in Brown and Walker\cite{bw} with the original $b$). Since the matrix is range symmetric,  the GMRES, RR-GMRES, and stabilized GMRES can be directly applied to $Ax=b$ (See \cite{bw} Theorem 2.4, \cite{hayami2011g} Theorem 2.7, and \cite{cr} Theorem 3.2.) as shown in Table \ref{tbbw42}. The stabilized GMRES gave a relative residual norm 1.94$\times 10^{-11}$ for bw42 at the 219th iteration. The proposed method can be considered as a way of making the GMRES stable for highly ill-conditioned inconsistent problems.
\begin{figure}
  \centering
  \begin{tabular}{c}
  \includegraphics[width=11cm]{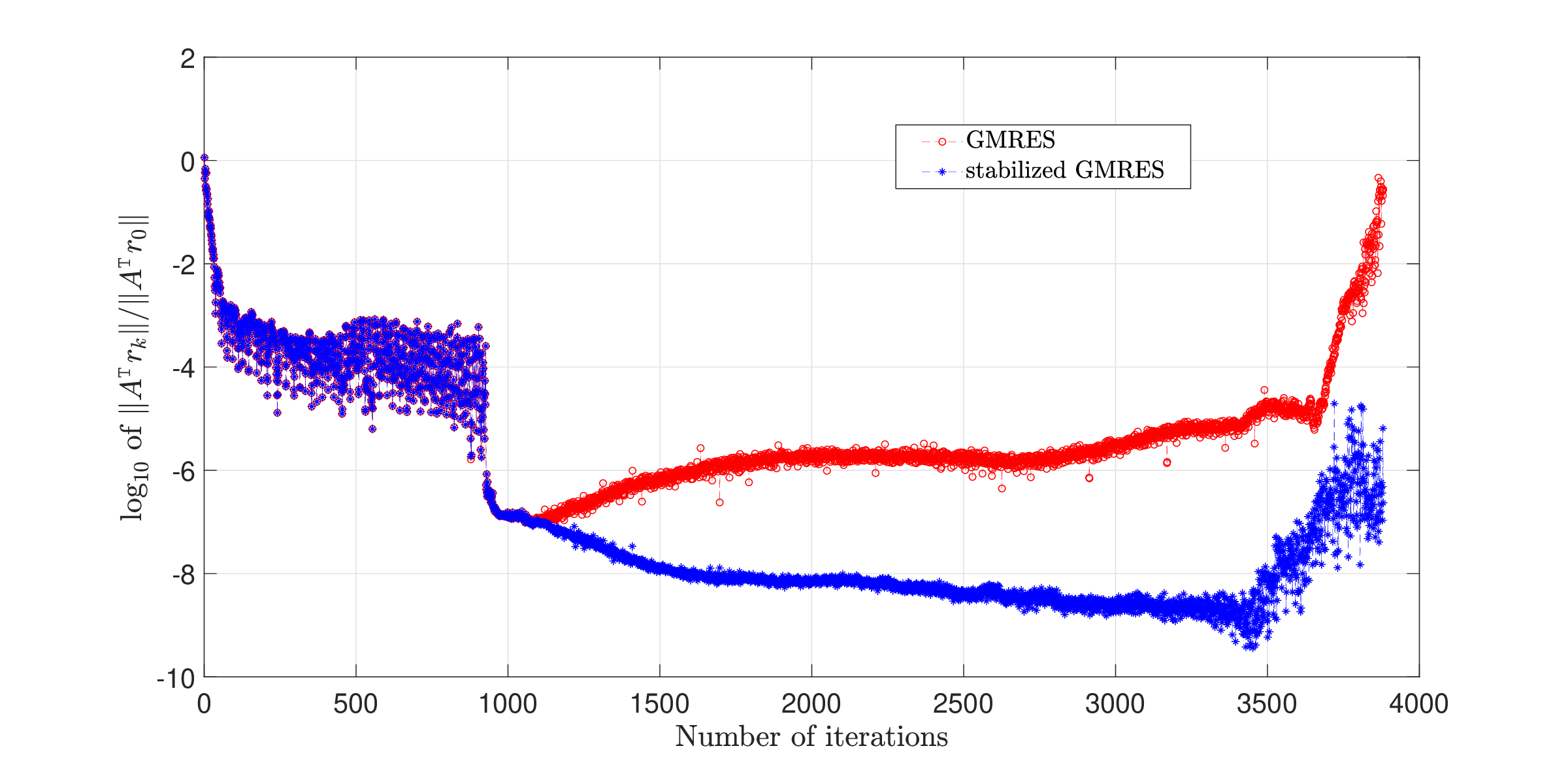}
\end{tabular}
  \caption{Comparison of GMRES with stabilized GMRES for freeFlyingRobot$\_$7.}
  \label{lwp12}
 \end{figure}
\par
\begin{figure}
  \centering
  \begin{tabular}{c}
  \includegraphics[width=11cm]{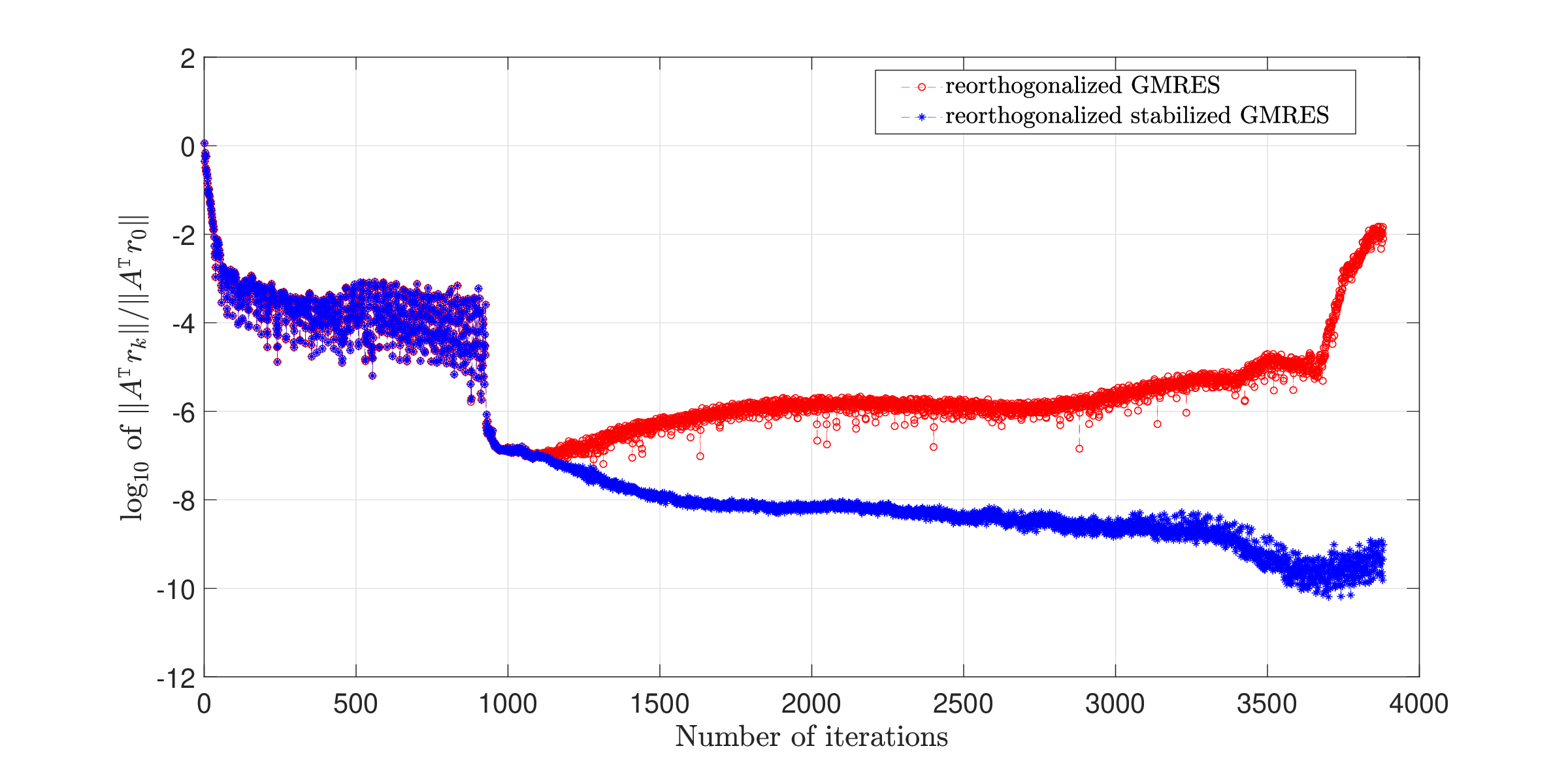}
\end{tabular}
  \caption{Comparison of reorthogonalized GMRES with reorthogonalized stabilized GMRES for freeFlyingRobot$\_$7.}
  \label{lwp11}
 \end{figure}

Figure \ref{lwp12} shows comparison of GMRES with stabilized GMRES for a symmetric matrix (, which is range symmetric), freeFlyingRobot$\_$7 which contains a cluster of tiny singular values which gradually decrease to zero. The stabilized GMRES converged to $3.65\times 10^{-10}$ at 3,452 iterations, better than GMRES. But the relative residual increased after the 3,452 iterations. Hence, we adopted a reorthogonalization strategy which performs the modified Gram-Schmidt orthogonalization process once more. We replaced line 4-6 of Algorithm \ref{AL1} (GMRES version) by Algorithm \ref{AL3} to reorthogonalize GMRES and the stabilized GMRES. As in Figure \ref{lwp11}, after reorthogonalization, the stabilized GMRES became more stabilized and converged to a relative residual of $6.45\times 10^{-11}$ at 3,701 iterations.
\begin{algorithm}
\caption{reorthogonalized modified Gram-Schmidt}
\label{AL3}
\begin{algorithmic}[1]
\For{$i=1,2$}
\For{$j=1,2,\dots,k$}
\State $h_{j, k}=w_k^{{\mathsf{T}}}v_j$, \quad  $w_k=w_k-h_{j, k}v_j$
\EndFor
\EndFor
\end{algorithmic}
\end{algorithm}
\par
\section{Concluding Remarks}
We proposed a stabilized AB-GMRES method for ill-conditioned underdetermined and inconsistent least squares problems. It shifts upwards the tiny singular values of the upper triangular matrix appearing in AB-GMRES, making the process more stable, giving better convergence, and more accurate solutions compared to AB-GMRES. We have also given a theoretical analysis to explain why the proposed method works. The method is also effective for making GMRES stable for range-symmetric inconsistent least squares problems with severely ill-conditioned square coefficient matrices.

\section*{Acknowledgments}
We would like to thank Dr. Hiroshi Murakami and Professor Lothar Reichel for valuable comments.\par
Ken Hayami was supported by JSPS KAKENHI Grant Number JP15K04768.\par
Keiichi Morikuni was supported by JSPS KAKENHI Grant Numbers JP16K17639 and JP20K14356 and Hattori Hokokai Foundation.\par
Jun-Feng Yin was supported by the National Natural Science Foundation of China (No. $11971354$).

\appendix
\section{Proof of statement in section 2.3}\label{ap1}
\begin{lem}\label{lemma1}
Assume $\rm{\mathcal{N}}$$(\hat{A})$ $\cap$ $\rm{\mathcal{R}}$$(\hat{\it{A}})=\{0\}$, and  $\rm{grade}$$(\hat{A},b|_{\rm{\mathcal{R}(\hat{\it{A}})}})=k$. Then, \\ $\mathcal{K}_{k+1}(\hat{\it{A}},b|_{\rm{\mathcal{R}(\hat{\it{A}})}})$ = $\hat{A}\mathcal{K}_{k}(\hat{\it{A}},b|_{\rm{\mathcal{R}(\hat{\it{A}})}})$ holds.
\begin{proof}
\rm {Note that}
\begin{align*}
\hat{\it{A}}\mathcal{K}_{k}(\hat{\it{A}},b|_{\rm{\mathcal{R}(\hat{\it{A}})}})&=\rm{span}\it\{\hat{\it{A}}b|_{\rm{\mathcal{R}(\hat{\it{A}})}}, \hat{\it{A}}^{\rm 2\it}b|_{\rm{\mathcal{R}(\hat{\it{A}})}}, \cdots,\hat{\it{A}}^k b|_{\rm{\mathcal{R}(\hat{\it{A}})}}\}\\
&\subseteqq \rm{span}\it\{b|_{\rm{\mathcal{R}(\hat{\it{A}})}},\hat{\it{A}}b|_{\rm{\mathcal{R}(\hat{\it{A}})}},\cdots,\hat{A}^k b|_{\rm{\mathcal{R}(\hat{\it{A}})}}\}=\mathcal{K}_{k+1}(\hat{\it{A}},b|_{\rm{\mathcal{R}(\hat{\it{A}})}}).
\end{align*}

\noindent grade$(\hat{\it{A}},b|_{\rm{\mathcal{R}(\hat{\it{A}})}})=k$ implies that

\begin{equation*}
\mathcal{K}_{k+1}(\hat{\it{A}},b|_{\rm{\mathcal{R}(\hat{\it{A}})}})=\mathcal{K}_{k}(\hat{\it{A}},b|_{\rm{\mathcal{R}(\hat{\it{A}})}})= \rm{span}\it\{b|_{\rm{\mathcal{R}(\hat{\it{A}})}},\hat{\it{A}}b|_{\rm{\mathcal{R}(\hat{\it{A}})}},\cdots,\hat{\it{A}}^{k-1} b|_{\rm{\mathcal{R}(\hat{\it{A}})}}\}.
\end{equation*}

\noindent Hence,
\begin{equation*}
\hat{\it{A}}^k b|_{\rm{\mathcal{R}(\hat{\it{A}})}}=c_0b|_{\rm{\mathcal{R}(\hat{\it{A}})}}+c_1\hat{\it{A}}b|_{\rm{\mathcal{R}(\hat{\it{A}})}}+\cdots+c_{k-1}\hat{\it{A}}^{k-1} b|_{\rm{\mathcal{R}(\hat{\it{A}})}},\quad c_i\in \mathbb{R}, i=0, 1, 2, \cdots, k-1.
\end{equation*}

\noindent If $c_0=0,$
\begin{equation*}
\hat{\it{A}}^k b|_{\rm{\mathcal{R}(\hat{\it{A}})}}=c_1\hat{\it{A}}b|_{\rm{\mathcal{R}(\hat{\it{A}})}}+c_2\hat{\it{A}}^2b|_{\rm{\mathcal{R}(\hat{\it{A}})}}+\cdots+c_{k-1}\hat{\it{A}}^{k-1} b|_{\rm{\mathcal{R}(\hat{\it{A}})}}.
\end{equation*}

\noindent Hence,
\begin{multline*}
c_1\hat{\it{A}}b|_{\rm{\mathcal{R}(\hat{\it{A}})}}+c_2\hat{\it{A}}^2b|_{\rm{\mathcal{R}(\hat{\it{A}})}}+\cdots+c_{k-1}\hat{\it{A}}^{k-1} b|_{\rm{\mathcal{R}(\hat{\it{A}})}}-\hat{\it{A}}^k b|_{\rm{\mathcal{R}(\hat{\it{A}})}}\\=\hat{\it{A}}(c_1b|_{\rm{\mathcal{R}(\hat{\it{A}})}}+\cdots+c_{k-1}\hat{\it{A}}^{k-2} b|_{\rm{\mathcal{R}(\hat{\it{A}})}}-\hat{\it{A}}^{k-1} b|_{\rm{\mathcal{R}(\hat{\it{A}})}})=0.
\end{multline*}

\noindent Hence,
\begin{equation*}
c_1b|_{\rm{\mathcal{R}(\hat{\it{A}})}}+c_2\hat{\it{A}}^2b|_{\rm{\mathcal{R}(\hat{\it{A}})}}+\cdots+c_{k-1}\hat{\it{A}}^{k-2} b|_{\rm{\mathcal{R}(\hat{\it{A}})}}-\hat{\it{A}}^{k-1} b|_{\rm{\mathcal{R}(\hat{\it{A}})}}\in \rm{\mathcal{N}}(\hat{\it{A}})\cap \rm{\mathcal{R}}(\hat{\it{A}}) = \{0\}.
\end{equation*}



\noindent which implies

\begin{equation*}
\hat{\it{A}}^{k-1} b|_{\rm{\mathcal{R}(\hat{\it{A}})}}=c_1b|_{\rm{\mathcal{R}(\hat{\it{A}})}}+c_2\hat{\it{A}}b|_{\rm{\mathcal{R}(\hat{\it{A}})}}+\cdots+c_{k-1}\hat{\it{A}}^{k-2} b|_{\rm{\mathcal{R}(\hat{\it{A}})}}.
\end{equation*}

\noindent Thus,
\begin{equation*}
\mathcal{K}_{k}(\hat{A},b|_{\rm{\mathcal{R}(\hat{\it{A}})}})=\mathcal{K}_{k-1}(\hat{\it{A}},b|_{\rm{\mathcal{R}(\hat{\it{A}})}}),
\end{equation*}
which contradicts with  grade$(\hat{A},b|_{\rm{\mathcal{R}(\hat{\it{A}})}})=k$. Hence, $c_0\neq 0$, and

\begin{equation*}
b|_{\rm{\mathcal{R}(\hat{\it{A}})}}=d_1\hat{\it{A}}b|_{\rm{\mathcal{R}(\hat{\it{A}})}}+d_2\hat{\it{A}}^2b|_{\rm{\mathcal{R}(\hat{\it{A}})}}+\cdots+d_{k-1}\hat{\it{A}}^{k-1} b|_{\rm{\mathcal{R}(\hat{\it{A}})}}+d_k\hat{\it{A}}^k b|_{\rm{\mathcal{R}(\hat{\it{A}})}}.
\end{equation*}

\noindent Hence,
\begin{multline*}
\mathcal{K}_{k+1}(\hat{A},b|_{\rm{\mathcal{R}(\hat{\it{A}})}})=\rm{span}\it\{b|_{\rm{\mathcal{R}(\hat{\it{A}})}},\hat{A}b|_{\rm{\mathcal{R}(\hat{\it{A}})}},\cdots,\hat{A}^k b|_{\rm{\mathcal{R}(\hat{\it{A}})}}\}\\ \subseteqq\rm{span}\it\{\hat{\it{A}}b|_{\rm{\mathcal{R}(\hat{\it{A}})}}, \hat{\it{A}}^{\rm 2\it}b|_{\rm{\mathcal{R}(\hat{\it{A}})}}, \cdots,\hat{A}^k b|_{\rm{\mathcal{R}(\hat{\it{A}})}}\}=\hat{A}\mathcal{K}_{k}(\hat{A},b|_{\rm{\mathcal{R}(\hat{\it{A}})}}).
\end{multline*}
\noindent Thus,
\begin{equation*}
\mathcal{K}_{k+1}(\hat{A},b|_{\rm{\mathcal{R}(\hat{\it{A}})}})=\hat{A}\mathcal{K}_{k}(\hat{A},b|_{\rm{\mathcal{R}(\hat{\it{A}})}}).
\end{equation*}
\end{proof}
\end{lem}
\begin{cor}\label{corollary1}
Assume $\rm{\mathcal{N}}$$(\hat{A})$ $=$ $\rm{\mathcal{N}}$$(\hat{\it{A}}^{\mathsf{T}})$, and  $\rm{grade}$$(\hat{A},b|_{\rm{\mathcal{R}(\hat{\it{A}})}})=k.$ Then,\\ $\mathcal{K}_{k+1}(\hat{A},b|_{\rm{\mathcal{R}(\hat{\it{A}})}})=\hat{A}\mathcal{K}_{k}(\hat{A},b|_{\rm{\mathcal{R}(\hat{\it{A}})}})$ holds.
\begin{proof}
 $\rm{\mathcal{N}}$$(\hat{A})$ $=$ $\rm{\mathcal{N}}$$(\hat{\it{A}}^{\mathsf{T}})$ \rm{implies that}
\begin{equation*}
\rm{\mathcal{N}}(\hat{\it{A}})\cap \rm{\mathcal{R}}(\hat{\it{A}})=\rm{\mathcal{N}}(\hat{\it{A}}^{\mathsf{T}})\cap \rm{\mathcal{R}}(\hat{\it{A}})=\rm{\mathcal{R}}(\hat{\it{A}})^{\bot}\cap \rm{\mathcal{R}}(\hat{\it{A}})=\{0\}.
\end{equation*}
Hence, from Lemma \ref{lemma1}, Corollary \ref{corollary1} holds.
\end{proof}
\end{cor}
\section{Proof of statement in section 2.3}\label{ap2}
\begin{lem}\label{lemma2} Assume $\rm{\mathcal{N}}$$(\hat{A})\cap \rm{\mathcal{R}}(\hat{\it{A}})=\{0\}$,
 $\rm{grade}$$(\hat{A},b|_{\rm{\mathcal{R}(\hat{\it{A}})}})=k$, and  $b\notin \mathcal{R}(\hat{A})$. Then, $\rm{dim}$$(\mathcal{K}_{k+1}(\hat{A},b))=k+1$ holds.\\
\begin{proof}
\rm {Let} $c_0, c_1, \dots, c_k\in \mathbb{R}$ satisfy
\begin{equation*}
c_0b+c_1\hat{\it{A}}b+\cdots+c_{k}\hat{\it{A}}^{k} b=0.
\end{equation*}
\noindent Since $\rm{\mathcal{N}}$$(\hat{A})$ $\cap$ $\rm{\mathcal{R}}$$(\hat{\it{A}})=\{0\}$,
\begin{equation*}
b=b|_{\rm{\mathcal{R}}(\hat{\it{A}})}\oplus b|_{\rm{\mathcal{N}}(\hat{\it{A}})},
\end{equation*}
where $b|_{\rm{\mathcal{N}}(\hat{\it{A}})}$ denotes the orthogonal projection of $b$ onto $\mathcal{N}(\hat{\it{A}}).$ Hence,
\begin{equation*}
c_0b|_{\rm{\mathcal{N}(\hat{\it{A}})}}+c_0b|_{\rm{\mathcal{R}(\hat{\it{A}})}}+c_1\hat{\it{A}}b|_{\rm{\mathcal{R}(\hat{\it{A}})}}+\cdots+c_{k}\hat{\it{A}}^{k} b|_{\rm{\mathcal{R}(\hat{\it{A}})}}=0.
\end{equation*}

\noindent If $c_0\neq 0$
\begin{equation*}
b|_{\rm{\mathcal{N}(\hat{\it{A}})}}=-b|_{\rm{\mathcal{R}(\hat{\it{A}})}}-\frac{c_1}{c_0}\hat{A}b|_{\rm{\mathcal{R}(\hat{\it{A}})}}-\cdots-\frac{c_{k}}{c_0}\hat{A}^{k} b|_{\rm{\mathcal{R}(\hat{\it{A}})}}\in \rm{\mathcal{R}(\hat{\it{A}})}.
\end{equation*}

\noindent Hence,
\begin{equation*}
b|_{\rm{\mathcal{N}(\hat{\it{A}})}}\in \rm{\mathcal{N}}(\hat{\it{A}})\cap \rm{\mathcal{R}}(\hat{\it{A}})=\{0\}.
\end{equation*}
Thus, $b|_{\rm{\mathcal{N}(\hat{\it{A}})}}=0,$ which contradicts $b\notin \mathcal{R}(\hat{\it{A}})$. Hence, we have $c_0= 0$, and
\begin{equation*}
c_1\hat{A}b+c_2\hat{A}^2b+\cdots+c_{k}\hat{A}^{k} b=c_1\hat{A}b|_{\rm{\mathcal{R}(\hat{\it{A}})}}+c_2\hat{A}^2b|_{\rm{\mathcal{R}(\hat{\it{A}})}}+\cdots+c_{k}\hat{A}^{k} b|_{\rm{\mathcal{R}(\hat{\it{A}})}}=0.
\end{equation*}

\noindent But, since
\begin{multline*}
\rm dim( \rm{span}\it\{\hat{\it{A}}b|_{\rm{\mathcal{R}(\hat{\it{A}})}}, \hat{A}^{\rm 2\it}b|_{\rm{\mathcal{R}(\hat{\it{A}})}}, \cdots,\hat{A}^k b|_{\rm{\mathcal{R}(\hat{\it{A}})}}\})\\=\rm dim( \rm\hat{\it{A}}\, {span}\it\{b|_{\rm{\mathcal{R}(\hat{\it{A}})}},\hat{A}b|_{\rm{\mathcal{R}(\hat{\it{A}})}}\cdots,\hat{A}^{k-1} b|_{\rm{\mathcal{R}(\hat{\it{A}})}}\})=\rm dim(\it \hat{\it{A}}\mathcal{K}_{k}(\hat{A},b|_{\rm{\mathcal{R}(\hat{\it{A}})}}))=k
\end{multline*}

\noindent holds from Lemma \ref{lemma1}, we have $c_1=c_2=\cdots=c_k=0,$ which implies $ \rm dim(\it \mathcal{K}_{k+1}(\hat{\it{A}},b))=k+1$.
\end{proof}
\end{lem}
\begin{cor} Assume $\rm{\mathcal{N}}$$(\hat{\it{A}})=$ $\rm{\mathcal{N}}$$(\hat{A}^{\mathsf{T}})$,  $\rm{grade}$$(\hat{A},b|_{\rm{\mathcal{R}(\hat{\it{A}})}})=k$, and  $b\notin \mathcal{R}(\hat{A})$. \\Then, $\rm{dim}$$(\mathcal{K}_{k+1}(\hat{A},b))=k+1$ holds.\\
\begin{proof}
$\rm{\mathcal{N}}$$(\hat{\it{A}})=$ $\rm{\mathcal{N}}$$(\hat{A}^{\mathsf{T}})$ \rm{implies} $\rm{\mathcal{N}}$$(\hat{A})$ $\cap$ $\rm{\mathcal{R}}$$(\hat{\it{A}})=\{0\}$. Hence, the corollary follows from Lemma \ref{lemma2}.
\end{proof}
\end{cor}
\section{Proof of Theorem \ref{th55} in section 4.4}\label{ap3}

\noindent \textbf{Theorom 9}\it\quad
Let $R_k=(r_{pq})\in \mathbb{R}^{k\times k}$ be an upper-triangular matrix and
\begin{equation}
   R_{k+\rm{1}\it}=\left(
             \begin{array}{cc}
               R_k & d \\
              0^{\mathsf{T}} & r_{k+\rm{1}\it,k+\rm{1}\it} \\
             \end{array}
           \right)\in \mathbb{R}^{(k+\rm{1}\it)\times(k+\rm{1}\it)}.\end{equation}
 Assume $R_k$ is numerically nonsingular, and $R_k=\mathcal{O}(\rm{1}\it)$, $ R_k^{-\rm{1}\it}=\mathcal{O}(\rm{1}\it)$, $ M(R_k)^{-\rm{1}\it}=\mathcal{O}(\rm{1}\it)$, $ d=\mathbb{O}(\rm{1}\it)$ and $O(k)=O(k^2)=O(\rm{1}\it).$ Then, the following holds:
\begin{equation*}
\rm f\!l\it(R_{k+\rm{1}\it}^{\mathsf{T}}R_{k+\rm{1}\it})~is~numerically~nonsingular \it\quad\Longleftrightarrow\quad \rm f\!l \it (r^2_{k+\rm{1}\it,k+\rm{1}\it})> \rm f\!l \it(d^{\mathsf{T}}d)O(u).
\end{equation*}

\begin{proof}\par\rm
Note that
\begin{equation*}
 R_{k+\rm{1}\it}^{\mathsf{T}}R_{k+\rm{1}\it}=\left(
             \begin{array}{cc}
               R_k & 0 \\
               d^{\mathsf{T}} & r_{k+\rm{1}\it,k+\rm{1}\it} \\
             \end{array}
           \right)\left(
             \begin{array}{cc}
               R_k &   d \\
            0^{\mathsf{T}}  & r_{k+\rm{1}\it,k+\rm{1}\it} \\
             \end{array}
           \right)=\left(
             \begin{array}{cc}
               R_k^{\mathsf{T}}R_k & R_k^{\mathsf{T}}d \\
               d^{\mathsf{T}}R_k & d^{\mathsf{T}}d+r_{k+\rm{1}\it,k+\rm{1}\it}^2 \\
             \end{array}
           \right).
    \end{equation*}
\\ \hspace*{\fill} \\
Proof of ($\Rightarrow$)
\\ \hspace*{\fill} \\
Assume  fl$(r^2_{k+\rm{1}\it,k+\rm{1}\it})$ $\leq$ fl$(d^{\mathsf{T}}d)O(u).$ Then, since \\
\begin{align*}
\rm f\!l\it (d^{\mathsf{T}}d)&=d^{\mathsf{T}}d+O(ku)d^{\mathsf{T}}d=(\rm{1}\it+O(ku))d^{\mathsf{T}}d,\\
\rm f\!l\it(d^{\mathsf{T}}d+r^{\rm 2\it}_{k+\rm{1}\it,k+\rm{1}\it})&=(d^{\mathsf{T}}d+r^2_{k+\rm{1}\it,k+\rm{1}\it})(\rm{1}\it+O(ku))=d^{\mathsf{T}}d(\rm{1}\it+O(ku)),\\
R_k&=\mathcal{O}(\rm{1}\it),\quad \rm and\quad \it d=\mathbb{O}(\rm{1}\it),
\end{align*}
we have\\
\begin{align}\label{eqapc1}
\rm f\!l\it(R_{k+\rm{1}\it}^{\mathsf{T}}R_{k+\rm{1}\it})&=\left(
             \begin{array}{cc}
               R_k^{\mathsf{T}}R_k +O(ku)|R_k|^{\mathsf{T}}|R_k|& R_k^{\mathsf{T}}d+O(ku)|R_k|^{\mathsf{T}}|d| \\
               d^{\mathsf{T}}R_k+O(ku)|d|^{\mathsf{T}}|R_k| & d^{\mathsf{T}}d+O(ku)d^{\mathsf{T}}d\\
             \end{array}
           \right)\notag\\&=\left(
             \begin{array}{c}R_k^{\mathsf{T}}\\ d^{\mathsf{T}}\end{array}\right)\left(
             \begin{array}{cc}R_k & d\end{array}\right)+\mathcal{O}(ku).
\end{align}
Note
\begin{equation*}\left(
             \begin{array}{cc}R_k & d\end{array}\right)\left(
             \begin{array}{c}-R_k^{-\rm{1}\it}d\\ \rm{1}\it\end{array}\right)=-R_k R_k^{-\rm{1}\it}d+d=0,
\end{equation*}
since $R_k$ is nonsingular.\par
Hence,\\
\begin{equation*}\rm f\!l\it(\left(
             \begin{array}{cc}R_k & d\end{array}\right)\left(
             \begin{array}{c}-R_k^{-\rm{1}\it}d\\ \rm{1}\it\end{array}\right))=\rm f\!l\it \{R_k  \rm f\!l\it(-R_k^{-\rm{1}\it}d)+d\} = [\rm f\!l \it \{R_k\rm f\!l\it(-R_k^{-\rm{1}\it}d)\}+d]\{\rm{1}\it+O(u)\}.
\end{equation*}
Note here that
\begin{equation*}
\rm f\!l\it \{R_k\rm f\!l\it (-R_k^{-\rm{1}\it}d)\}=R_k\rm f\!l\it(-R_k^{-\rm{1}\it}d)+O(ku)|R_k||R_k^{-\rm{1}\it}d|,
\end{equation*}
and
\begin{equation}\label{eqapc2}
\rm f\!l\it (-R_k^{-\rm{1}\it}d)= -R_k^{-\rm{1}\it}d+O(k^{\rm 2\it}u)M(R_k)^{-\rm{1}\it}|d|
\end{equation}
from Theorem \ref{th333}. Hence,
\begin{equation*}
\rm f\!l\it (\left(
             \begin{array}{cc}R_k & d\end{array}\right)\left(
             \begin{array}{c}-R_k^{-\rm{1}\it}d\\ \rm{1}\it\end{array}\right)) =
 O(k^{\rm 2\it}u)R_sM(R_k)^{-\rm{1}\it}|d|+O(ku)|R_k||R_k^{-\rm{1}\it}d| =\mathbb{O}(k^{\rm 2\it}u),
\end{equation*}
since $R_k^{-\rm{1}\it}=\mathcal{O}(\rm{1}\it)$ and $M(R_k)^{-\rm{1}\it}=\mathcal{O}(\rm{1}\it).$\par
Then,
\begin{multline*}
\rm f\!l\it(R_{k+\rm{1}\it}^{\mathsf{T}}R_{k+\rm{1}\it}\left(
             \begin{array}{c}-R_k^{-\rm{1}\it}d\\ \rm{1}\it\end{array}\right))\\ = \rm f\!l\it(\{\left(
             \begin{array}{c}R_k^{\mathsf{T}}\\ d^{\mathsf{T}}\end{array}\right)\left(
             \begin{array}{cc}R_k & d\end{array}\right)+\mathcal{O}(ku)\}\left(
             \begin{array}{c}-R_k^{-\rm{1}\it}d+\mathcal{O}(k^2u)M(R_k)^{-\rm{1}\it}|d|\\ \rm{1}\it\end{array}\right))=\mathbb{O}(k^{\rm 2\it}u)=\mathbb{O}(u),
\end{multline*}
since (\ref{eqapc1}), (\ref{eqapc2}), and $O(k^2)=O(\rm{1}\it).$ Since $\left(\begin{array}{c}-R_k^{-\rm{1}\it}d\\ \rm{1}\it\end{array}\right)=\mathbb{O}(\rm{1}\it),$ $R_{k+\rm{1}\it}^{\mathsf{T}}R_{k+\rm{1}\it}$ is numerically singular. By contraposition, ($\Rightarrow$)~holds.
\\ \hspace*{\fill} \\

\noindent Proof of ($\Leftarrow$)
\\ \hspace*{\fill} \\
Assume  $R_{k+\rm{1}\it}^{\mathsf{T}}R_{k+\rm{1}\it}$ is not numerically nonsingular. Then, there exists a vector $\left(
                                                                                                        \begin{array}{c}
                                                                                                          z \\
                                                                                                          w \\
                                                                                                        \end{array}
                                                                                                      \right)\in \mathbb{R}^{k+\rm{1}\it}$ such that $\left|\left(
                                                                                                        \begin{array}{c}
                                                                                                          z \\
                                                                                                          w \\
                                                                                                        \end{array}
                                                                                                      \right)\right|>\mathbb{O}(u),$
and \begin{align*}
\rm f\!l\it \{R_{k+\rm{1}\it}^{\mathsf{T}}R_{k+\rm{1}\it}\left(
                                                                                                        \begin{array}{c}
                                                                                                          z \\
                                                                                                          w \\
                                                                                                        \end{array}
                                                                                                      \right)\}= R_{k+\rm{1}\it}^{\mathsf{T}}\left(R_{k+\rm{1}\it}\left(
                                                                                                        \begin{array}{c}
                                                                                                          z \\
                                                                                                          w \\
                                                                                                        \end{array}
                                                                                                      \right)+|R_{k+\rm{1}\it}|\left|\left(
                                                                                                        \begin{array}{c}
                                                                                                          z \\
                                                                                                          w \\
                                                                                                        \end{array}
                                                                                                      \right)\right|O((k+\rm{1}\it)u)\right)+\\ \left|R_{k+\rm{1}\it}^{\mathsf{T}}\right|\left|R_{k+\rm{1}\it}\left(
                                                                                                        \begin{array}{c}
                                                                                                          z \\
                                                                                                          w \\
                                                                                                        \end{array}
                                                                                                      \right)+|R_{k+\rm{1}\it}|\left|\left(
                                                                                                        \begin{array}{c}
                                                                                                          z \\
                                                                                                          w \\
                                                                                                        \end{array}
                                                                                                      \right)\right|O((k+\rm{1}\it)u)\right|O((k+\rm{1}\it)u)=\mathbb{O}(u)
\end{align*}
 assuming $O(k+\rm{1}\it)=O(\rm{1}\it).$\\ Hence,
\begin{equation*} \rm f\!l\it\{R_{k+\rm{1}\it}^{\mathsf{T}}R_{k+\rm{1}\it}\left(
                                                                                                        \begin{array}{c}
                                                                                                          z \\
                                                                                                          w \\
                                                                                                        \end{array}
                                                                                                      \right)\}=\left(
             \begin{array}{cc}
               R_k^{\mathsf{T}}R_k & R_k^{\mathsf{T}}d \\
               d^{\mathsf{T}}R_k & d^{\mathsf{T}}d+r_{k+\rm{1}\it,k+\rm{1}\it}^2 \\
             \end{array}
           \right)\left(
                                                                                                        \begin{array}{c}
                                                                                                          z \\
                                                                                                          w \\
                                                                                                        \end{array}
                                                                                                      \right)+\mathbb{O}(u)=\mathbb{O}(u).\end{equation*}
Thus,
\begin{equation}\label{c1}
 R_k^{\mathsf{T}}R_sz+wR_k^{\mathsf{T}}d=\mathbb{O}(u),
\end{equation}
\begin{equation}\label{c2}
d^{\mathsf{T}}R_sz+(d^{\mathsf{T}}d+r_{k+\rm{1}\it,k+\rm{1}\it}^2)w=\mathbb{O}(u).
\end{equation}
(\ref{c1}) can be expressed as $R_k^{\mathsf{T}}(R_sz+wd)=\mathbb{O}(u).$
From Lemma \ref{lemma3}, $R_k^\mathsf{T}$ is numerically nonsingular, so that
\begin{equation}\label{c3}
R_sz+wd=\mathbb{O}(u).
\end{equation}
Hence, from (\ref{c2}), $d^{\mathsf{T}}R_sz+w(d^{\mathsf{T}}d+r_{k+\rm{1}\it,k+\rm{1}\it}^2)=d^{\mathsf{T}}(R_sz+wd)+wr_{k+\rm{1}\it,k+\rm{1}\it}^2=O(u).$
Thus, $wr_{k+\rm{1}\it,k+\rm{1}\it}^2=O(u)$. If $w=O(u),$ $R_sz=\mathbb{O}(u)$ from (\ref{c3}).
Since $R_k$ is numerically nonsingular, $z=\mathbb{O}(u),$
which contradicts with the assumption.\par
Hence, $|w|>O(u),$ so that $r_{k+\rm{1}\it,k+\rm{1}\it}^2=O(u),$ which gives
\begin{equation*}
\rm  f\!l\it(r_{k+\rm{1}\it,k+\rm{1}\it}^{\rm 2\it}) = O(u)\leq \rm f\!l\it( d^{\mathsf{T}}d)O(u).\end{equation*}
\end{proof}
\rm
\begin{lem}\label{lemma3}\it
Let $n=O(\rm{1}\it)$. If $A\in \mathbb{R}^{n\times n}$ is numerically nonsingular, and $A^{-\rm{1}\it}=\mathcal{O}(\rm{1}\it)$, then $A^{\mathsf{T}}$ is numerically nonsingular.
\begin{proof}
\rm If
\begin{equation*}
\rm f\!l\it (A^{\mathsf{T}}x) = A^{\mathsf{T}}x+\mathbb{O}(nu)|A^{\mathsf{T}}||x|=\mathbb{O}(nu),
\end{equation*}
then
\begin{equation*}
\rm f\!l\it(x^{\mathsf{T}}A) = x^{\mathsf{T}}A+\mathbb{O}^{\mathsf{T}}(nu)=\mathbb{O}^{\mathsf{T}}(nu).
\end{equation*}
Thus,
\begin{equation*}
\rm f\!l(\it x^{\mathsf{T}}Ay) = \rm f\!l\it (x^{\mathsf{T}}A)y+O(nu)|\rm f\!l\it(x^{\mathsf{T}}A)||y|=O(nu)
\end{equation*}
holds for all $y=\mathbb{O}(\rm{1}\it).$\par
For arbitrary $z=\mathbb{O}(\rm{1}\it)\in \mathbb{R}^{n},$ let \begin{equation*}y=A^{-\rm{1}\it}z=\mathbb{O}(\rm{1}\it).\end{equation*}
Then,
\begin{equation*}
\rm f\!l(\it Ay) = Ay+O(nu)|A||y|=z+O(nu)|A||y|.
\end{equation*}
Hence,
\begin{equation*}
z=\rm f\!l\it(Ay) + O(nu)|A||y| = \rm f\!l\it(Ay) + \mathbb{O}(nu).
\end{equation*}
Thus, we have
\begin{equation*}
\rm f\!l(\it x^{\mathsf{T}}z) = x^{\mathsf{T}}z + O(nu)|x|^{\mathsf{T}}|z| =\rm f\!l\it(x^{\mathsf{T}}Ay) + O(nu)= O(nu)
\end{equation*}	
for arbitrary $z=\mathbb{O}(\rm{1}\it)\in \mathbb{R}^n.$ Hence, $x=\mathbb{O}(u),$ so that $A^{\mathsf{T}}$ is numerically nonsingular.
\end{proof}
\end{lem}
\section{Proof of Theorem \ref{thmeighth} in section 4.5}\label{ap4}
\begin{proof}\rm
Let the singular value decomposition of $R_k$ be given by $R_k=U\Sigma V^{\mathsf{T}}\in \mathbb{R}^{k\times k},$ where $U, V$ are orthogonal matrices and $\Sigma=\rm diag (\sigma_{\rm 1}, \sigma_2, \dots,\it \sigma_k).$ Let $I_k\in \mathbb{R}^{k\times k}$ be the identity matrix. Then, we have $R_k'=\left( \begin{array}{c}
R_k\\
\sqrt{\lambda}I_k\\
   \end{array}\right)=U'\Sigma'V^{\mathsf{T}}$, where $U'=\left( \begin{array}{cc}
U & 0\\
0 & V
   \end{array}\right)$ and $\Sigma'=\left( \begin{array}{c}
\Sigma\\
\sqrt{\lambda}I_k\\
   \end{array}\right).$
Since $\Sigma'^{\mathsf{T}}\Sigma'=\Sigma^2+\lambda I_k=\rm diag (\sigma_1^2+\lambda, \sigma_2^2+\lambda, \dots, \sigma_{\it k}^{\rm 2}+\lambda),$ the singular values of $\left( \begin{array}{c}
R_k\\
\sqrt{\lambda}I_k\\
   \end{array}\right)$ are  $\sqrt{\sigma_1^2+\lambda}\geq\sqrt{\sigma_2^2+\lambda}\geq\cdots\geq\sqrt{\sigma_k^2+\lambda}.$
\end{proof}

\end{document}